# Euler: Genius Blind Astronomer Mathematician

Dora E. Musielak
University of Texas at Arlington
dmusielak@uta.edu

**Résumé**. Leonhard Euler, the most prolific mathematician in history, contributed to advance a wide spectrum of topics in celestial mechanics. At the St. Petersburg Observatory, Euler observed sunspots and tracked the movements of the Moon. Combining astronomical observations with his own mathematical genius, he determined the orbits of planets and comets. Euler laid the foundations of the methods of planetary perturbations and solved many of the Newtonian mechanics problems of the eighteenth century that are relevant today. In his study of the three-body problem, Euler discovered two of five equilibrium points so-called the Lagrangian points. His pioneering work in astronomy was recognized with six of the twelve prizes he won from the Paris Academy of Sciences. In this article, we review some of Euler's most interesting work in astronomy.

**Résumé**. Leonard Euler, el matemático más prolífico de la historia, contribuyó al avance de una amplia gama de temas en la mecánica celeste. En el Observatorio de San Petersburgo, Euler observó las manchas solares e hizo observaciones de los movimientos de la Luna. Combinando las observaciones astronómicas con su propio genio matemático, él determinó las órbitas de los planetas y los cometas. Euler sentó las bases de los métodos de perturbaciones planetarias y resolvió muchos de los problemas de la mecánica newtoniana del siglo XVIII que aun hoy en día son relevantes. En su estudio sobre el problema de los tres cuerpos Euler descubrió dos de los puntos de equilibrio también llamados puntos de LaGrange. Su trabajo pionero en la astronomía fue reconocido con seis de los doce premios que ganó en la Academia de Ciencias de París. En este artículo se revisan algunos de los trabajos más interesantes de Euler en astronomía.

## 1. Introduction

Let us imagine Leonhard Euler at twenty-three, a passionate young man full of wonderment, observing the Moon and the stars through a telescope. He was a gifted mathematician who would go on to pioneer and develop the most beautiful and useful mathematical methods we use today. In 1730, before he lost his eyesight, Euler was an astronomer at the royal Observatory in Saint Petersburg, Russia.

For the first ten years of his career, Euler carried out astronomical observations as part of his research work at the St. Petersburg Academy of Sciences. While establishing his reputation as a mathematician, Euler observed sunspots and made direct observations of the movements of the Moon. In his first article related to astronomy, Euler introduced spherical triangles to obtain a method for determining the coordinates of a polar star.

Three years later, Euler dealt with the motion of planets and their orbits, in which he developed an iterative method of solution for Kepler's problem.

When his vision began to fail him, Euler stopped his daily routine at the observatory. However, he continued his research to derive the mathematical theories required to solve the most difficult and complex problems of celestial mechanics.[1] In the first years of his career, Euler created analytical mechanics, the refined and highly mathematical form of classical Newtonian mechanics[2] that describes the motion of macroscopic objects, including astronomical bodies, such as planets, stars, and galaxies. The new tools Euler developed served to solve some of the most pressing problems of astronomy of the eighteenth century.

With his analytical work, Euler also contributed to improving solar and planetary tables. After Kepler, these tables became Newtonian by incorporating the lunar and planetary perturbations implied by Newton's law of universal gravitation.[3] Essential to this transformation were the analytical results achieved by Euler, and after him Clairaut, d'Alembert, Lagrange, and Laplace. Analytical mechanics gave astronomy the mathematical foundation to elucidate the most difficult and complex problems of celestial mechanics that are relevant especially today.

The scope of Euler's contributions to astronomy embraced many areas: Astronomical Perturbation; Planetary Motion; Precession and Nutation; Eclipses and Parallax; Tides and Geophysics; and Optics, which included refinement of lenses for telescopes.[4] Euler's accomplishments include determining with great accuracy the orbits of comets and other celestial bodies, and calculating the parallax of the Sun. He also made direct observations of the Moon and Sun. Euler's passion for the stars never diminished, even when he lost vision in both eyes. Just before he died, the seventy-six blind mathematician pondered the discovery of a new planet and would have calculated its orbit if he had had enough time.

## 2. The First Years

Many biographies have been written about Euler, and many more are yet to be written. It is not my intention to attempt to compete with those works. I just wish to illuminate a few facts about Euler's life to help me chronicle his contributions to celestial mechanics.

### Beginning in Switzerland

Leonhard Euler was born in Basel on 15 April 1707. He was the son of a Protestant minister who had studied mathematics under Jakob Bernoulli, the first of the famous family of renowned mathematicians. Euler received basic instruction from his father who expected him to pursue a career in theology. At thirteen, Euler entered the University of Basel, where he attended the freshman course of Johann Bernoulli I (younger brother of Jakob), which included geometry and arithmetic.

---

[1] Celestial mechanics, named by Laplace in 1798 (*Mécanique Céleste*), is concerned with the perturbations on the motion of a body orbiting a larger one by a third body.
[2] The term "classical mechanics" was coined in the early twentieth century to describe the system of mathematical physics begun by Isaac Newton and improved by Euler, Lagrange, and others after them.
[3] Planetary Astronomy from the Renaissance to the Rise of Astrophysics, Eds. René Taton and Curtis Wilson, Cambridge University Press (August 25, 1995).
[4] Euler's memoirs were retrieved from The Euler Archive directed by Dominic Klyve (Central Washington University), Lee Stemkoski (Adelphi University), and Erik Tou (Pacific Lutheran University), and hosted by the Mathematical Association of America at http://eulerarchive.maa.org/.

During Saturday visits at his home, Johann Bernoulli became Euler's mentor, instructing and guiding the young boy outside the university's classroom. During that time Euler became friends with Johann's sons, Niklaus II and Daniel who were a few years older than him. This close contact with notable mathematicians (the great master and his sons), became very important to Euler's intellectual development. At sixteen, Euler obtained the academic degree of *magister* in philosophy. He then enrolled at the Theological Faculty.

Meanwhile in Russia, Peter the Great founded the Saint Petersburg Academy of Sciences, inspired by German mathematician and philosopher Gottfried Wilhelm von Leibniz; it was implemented by the Senate decree of February 8 (January 28 o.s.), 1724. The Russian Academy hired Daniel and Niklaus Bernoulli as professors, along with mathematician Jakob Hermann (a distant relative of Euler). The brothers left Basel with the promise of securing a position for Euler. The official invitation from Empress Catherina of Russia came trough Daniel Bernoulli at the beginning of the winter of 1726. While waiting to depart for Russia, the nineteen-year old Euler enrolled at the Medical Faculty in Basel to study medicine.[5]

At the same time, Euler applied for a professorship of physics at the University of Basel. To demonstrate his scientific acumen and research interests, he submitted a *Dissertatio physica de sono*, an essay in which the young man introduced his theory to explain the nature of the atmosphere. His premise was based on the theory of elasticity he had learned from his mentor, Johann Bernoulli. Euler also stated (without giving a proof) a formula for the speed of sound propagation, and obtained numerical values that are of the correct order of magnitude for air.[6] Despite the merits of his 16-page dissertation, Euler did not get the position. This turn of events sealed his future, for then Euler was free to accept the call from Russia.

## Voyage to Russia in 1727

In April 1727, a few days after Isaac Newton died, Leonhard Euler began his long journey to St. Petersburg to start his career at the Russian Academy of Sciences. Euler had just turned twenty years old.

From Basel to Saint Petersburg, Euler traveled by land and by water, navigating by boat along the Rhine River, and sailing on a ship through the Baltic Sea. The long journey took him through several major cities and ports of today's Germany, Poland, Estonia, and Russia. I do not know if Euler recorded the experiences of his voyage, but if he did it may have described a rather harrowing experience. His trip lasted several weeks, from early April to mid May 1727. As an old man, Euler dictated to his son, Johann Albrecht, a short autobiography where he recalled that first voyage:

" ... I set out from Basel right at the beginning of April, but arrived in Lubec at such an early time that no ship was ready to sail to Petersburg: I was forced therefore to take a ship going to Reval, and since the trip took almost four weeks, I soon found in Reval a ship to Stettin which brought me to Cronstadt. There I arrived on precisely the day when the death of the Empress Catherina I Alekseyevna became known ...”[7]

---

[5] Fellmann, E.A., *Leonhard Euler*, Birkhäuser Verlag, Switzerland (2007).p. 6

[6] Euler, L., *Dissertatio physica de sono* (Physical dissertation on sound). Published in Basel in 1727.

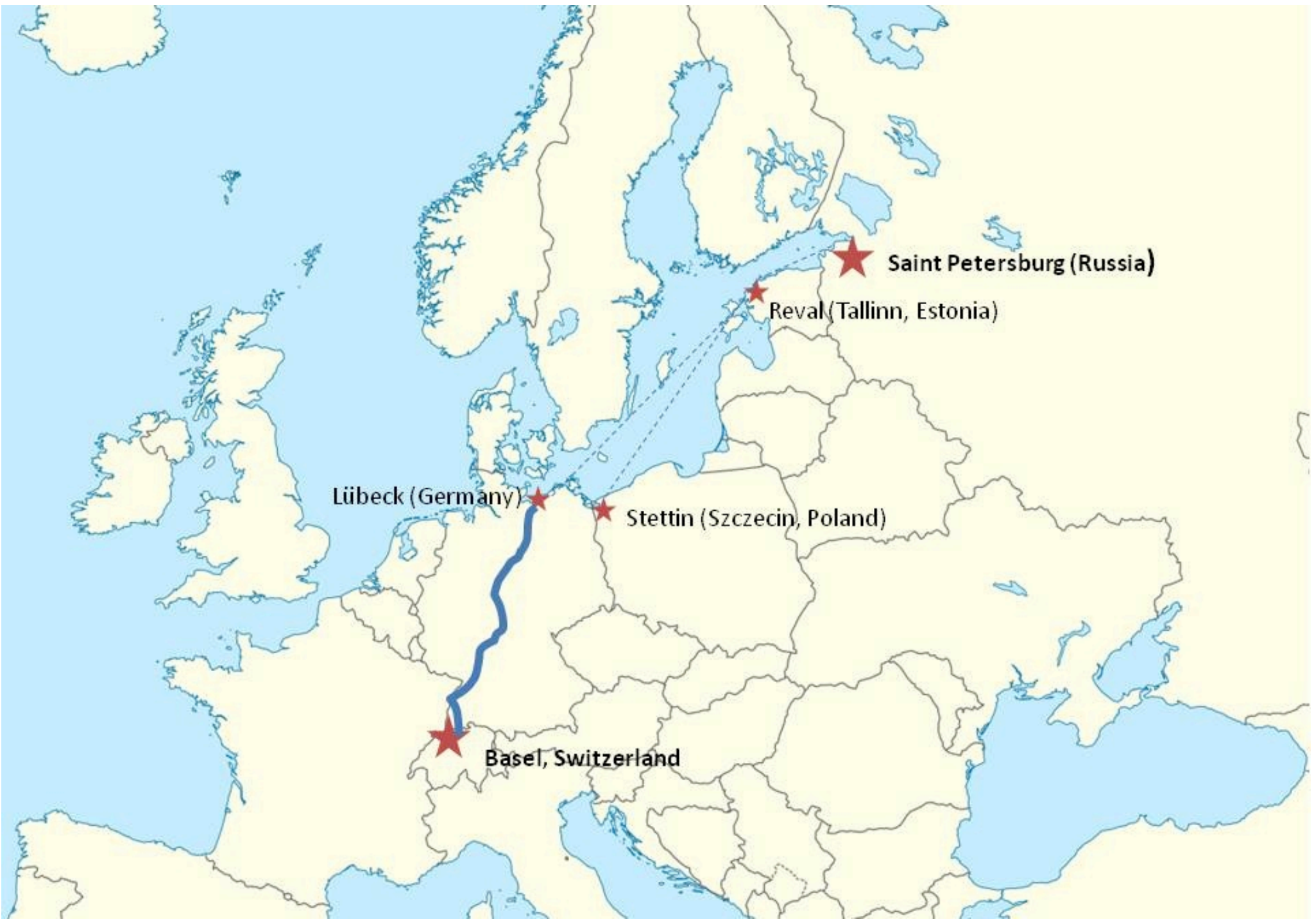

Modern map showing the path of Euler’s trip from Basel (Switzerland) to Saint Petersburg (Russia).

Perhaps Euler meant to say that from Lübec (a city port in today Germany) he sailed to Stettin (Szczecin, a port in today Poland). Then from Stettin he would take a ship to Reval (Tallinn, the capital city of Estonia), and from Reval to Cronstadt (Kronstadt), St. Petersburg’s main seaport, which is located on Kotlin Island, 30 kilometers (19 mi) west of St. Petersburg proper. Or maybe from Lübec he sailed straight to Reval, as it makes no sense that he would have returned to Stettin. (See the map indicating the path of Euler’s voyage.) Euler stepped on Russian soil on 17 May (6 May o.s.) 1727.

Travelling in the eighteenth century was rather difficult and strenuous. Did Euler walk some parts of his arduous journey? Or did he travel some tracks by wagon or carriage? The noble and the rich could travel in some comfort—in private, and in upholstered carriages accompanied by footmen and liveried coachmen of their own. But Euler was a common young man. The distance to be covered from Basel to Lübec is about 870 km. Today, the same distance can comfortably and easily be covered via the Intercity in about eight hours. At that time, it took weeks.

When Euler left Switzerland, Saint Petersburg was the imperial capital of the Russian Empire, a major European power ruled by tsar Peter the Great, Pyotr

[7] Translated from the German language by Emil A. Fellmann in “Leonhard Euler,” Birkhauser, 2007, p. 6.

Alexeyevich. Seen the strategic importance of the Baltic Sea, in 1703 Peter established his new capital at the mouth of the Neva River at the east end of the Gulf of Finland.

## St. Petersburg Academy and Its Observatory 1727-1738

The first years of his career at St. Petersburg must have been rewarding for Euler. His salary was 300 rubles and was given free lodging, firewood, and light. Among his mentors and friends he had Daniel Bernoulli, Christian Goldbach, Jakob Hermann,[8] and French astronomer Joseph-Nicolas Delisle. In 1730, Euler was professor of physics and, when Daniel returned to Basel, in 1733, Euler took over his professorship of mathematics. Euler must have been happy as well. On 7 January 1734 (27 December 1733 O.S.), he married Katharina Gsell. The young couple purchased a comfortable house on the banks of the Neva not far from the Academy. Their first child, Johann Albrecht, was born on 27 November 1734.

For almost ten years, Euler was among the astronomers taking measurements twice daily at the St. Petersburg Observatory. This is according to the observation records from 1725-1746 found in 1977 at the Leningrad archives of the Russian Academy of Sciences.[9] Those records reveal that Euler's entries were so detailed and numerous that some scholars deduce from them that Euler must have mastered the methods of astronomical observations.[10]

Founded in 1725, the royal observatory at St. Petersburg had a frontage of two hundred and twenty-five feet, and central towers one hundred and forty feet high.[11] The observatory was equipped with every variety of instruments known to the researchers of that time. French astronomer Joseph-Nicolas Delisle was recruited to create and run the school of astronomy. He arrived in Russia in 1726, after the death of the Russian czar Peter the Great, and became the Director of the St. Petersburg Observatory.

Upon his arrival, the observatory must have been captivating to the youthful Euler; perhaps it was the first he ever saw. He took the opportunity to learn to use the instruments of astronomy with Delisle.[12] It seems that sunspots were of much interest to Euler because his notes from this period contain enthusiastic comments on his observations.[13] Sunspots are dark spots on the Sun, some as large as 50,000 miles in diameter, which move across the surface, contracting and expanding as they go. These strange and powerful phenomena must have spiked Euler's curiosity. At that time, the nature of sunspots was speculative. Now we know that they are caused by the Sun's intense magnetic activity, which inhibits convection by an effect comparable to the eddy current brake, forming areas of reduced surface temperature.

---

[8] Swiss mathematician Jakob Hermann, a distant relative of Euler, worked on problems in classical mechanics. In 1716, he wrote *Phoronomia*, an early treatise on Mechanics. In St. Petersburg since 1724, Hermann returned to Basel in 1730 and died in 1733)

[9] Burckhardt, J.J., *Leonhard Euler 1707-1783*, Mathematics Magazine 56 (1983), 262-272. As it appeared in Sherlock Holmes in Babylon, Eds. M. Anderson, V. Katz, and R. Wilson, The Mathematical Association of America (2004).

[10] *Ibid*.

[11] Knowledge, An Illustrated Magazine of Science. Volume XVIII, January to December 1895. p. 87.

[12] Delisle published *Mémoires pour servir à l'histoire et au progrès de l'astronomie* (St. Petersburg, 1738), in which he gave the first method for determining the heliocentric coordinates of sunspots.

[13] Burckhardt, J.J., *Leonhard Euler 1707-1783*, Mathematics Magazine 56 (1983), 262-272. As it appeared in Sherlock Holmes in Babylon, Eds. M. Anderson, V. Katz, and R. Wilson, The Mathematical Association of America (2004).

Although sunspots had been observed since antiquity, it was the development of the telescope that began to provide a better understanding of their nature. Astronomers first observed sunspots telescopically in late 1610. Galileo, in particular, made them popular. Sunspots had some importance in the debate over the nature of the Solar System. They showed that the Sun rotated. Did Euler note the cyclic variation of the number of sunspots?

Applying findings from his twice-daily astronomical observations at the observatory, Euler wrote a number of memoirs of importance to contemporary astronomers. On 22 February 1732, Euler presented his first astronomy paper to the St. Petersburg Academy, in which he developed a mathematical method to calculate the coordinates of stars.[14] In this work, Euler used spherical trigonometry to introduce the spherical triangle theorem to prove that the cosine of a spherical triangle ABC is given by

$$\cos BC = \frac{\cos(AB + AC) + \cos(AB - AC)}{2} + \frac{\cos A \cos(AB - AC) - \cos A \cos(AB + AC)}{2}.$$

A spherical triangle is a triangular figure formed on the surface of a sphere by three great circular arcs intersecting pairwise in three vertices. Called an Euler triangle, the spherical triangle is the spherical analog of the planar triangle.

Euler stated the problem as follows (translated from Latin and illustrated by Ian Bruce):

> For a given fixed star successively observed in three places (Fig. 2) ABC, with the altitudes or the complements of these ZA, ZB, ZC, and with the elapsed times between the observations given, or from the angles to the pole P, APB and BPC, to find the elevation of the pole star or the complement PZ of this, and the declination of the star, or the complements of this, at either AP, BP or CP.

Euler meant that at some instant of time, a star's celestial coordinates are known by two angles: its right ascension or longitude in the celestial sphere, and the latitude or declination measured by the angles from the North Celestial Pole (pole star at P). With the spherical triangle theorem, Euler determined the elevation of the pole star and the declination of the star, knowing the star's elevation at different points in the sky. Astronomers refer to the angle from the horizon to a star as the angle of elevation.

At the end of the solution, Euler gave an example, which he "had computed from the altitude of the pole star in order that he might investigate whether or not the same numbers might arise from his method." Euler determined that it is possible to calculate the time for the motion of a star, and he concluded from it that "if the Sun itself is to take part in these observation, then the time of noon itself can be found." A way of defining solar noon is when the Sun's angle of elevation is greatest for a particular day. At such time, the length of shadows will be the shortest.

[14] Euler, L., *Solutio problematis astronomici ex datis tribus stellae fixae altitudinibus et temporum differentiis invenire elevationem poli et declinationem stellae* (Solution to problems of astronomy: given the altitudes and time differences for three fixed stars, to find the elevation of the pole and the declination of the star). February 22, 1732. E14. Translated and annotated by Ian Bruce.

Fig. 2

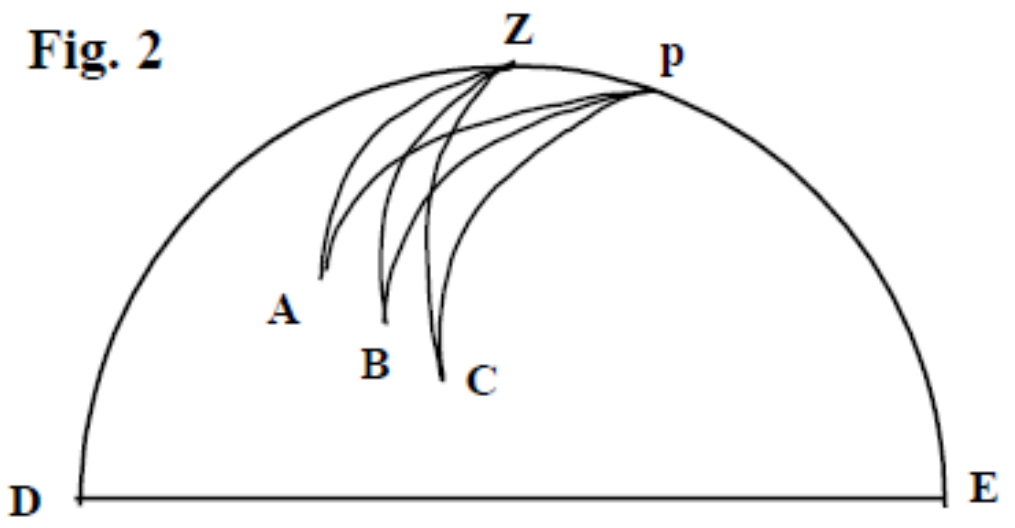

We should note that, in the 1700s, the investigation of solar noon was a very important task for astronomers. Solar noon is the moment when the Sun transits the celestial meridian—approximately the time when it is highest above the horizon on that day. Solar noon is not necessarily the same time as clock noon. For people, clock noon is a concept based on the notion that the interval from midday today to midday tomorrow ought to be the same every day of the year. Before accurate mechanical clocks were perfected, the concept of midday (solar noon) was derived from direct observation of the Sun's apparent movement across the sky. The true noon is the time or point in the Sun's path at which the Sun is on the local meridian. It is not surprising that the young Euler would have tested his mathematical skills formulating a simple theorem to calculate the precise instant of true noon.

On 28 February 1735, Euler presented to the St. Petersburg Academy a method for computing the equation of a meridian.[15] This equation was needed for the calculation of the meridian arc for the determination of the figure of the Earth. In geodesy, a meridian arc measurement is a highly accurate determination of the distance between two points with the same longitude. Two or more such determinations at different locations then specify the shape of the reference ellipsoid which best approximates the shape of the geoid.[16] The earliest determinations of the size of a spherical Earth required a single arc.

In 1687, Newton asserted that the Earth was an oblate spheroid, and scientists after him were disputing or attempting to prove it. Euler alone found a simple method of computing tables for the meridional equation of the Sun that was important to confirm Newton's prediction that the Earth was slightly flattened at the poles. It was not until the next year when the French expedition to Lapland led by Pierre-Louis de Maupertuis (accompanied by mathematician Alexis Clairaut and other scientists) that Newton's theory was fully confirmed. In 1738, Maupertuis published *Sur la figure de la terre*, the memoir where he described the expedition and the meridian measurements, giving conclusive data to support Newton's assertion.

[15] Euler, L., *Methodus computandi aequationem meridiei* (E50), February 28, 1735. Published in Commentarii academiae scientiarum Petropolitanae 8, 1741, pp. 48-65

[16] The geoid is the shape that the surface of the oceans would take under the influence of Earth's gravitation and rotation alone, in the absence of other influences such as winds and tides.

At the end of 1735, Euler published his work on the motion of planets,[17] where he calculated their orbits. There we find the familiar orbit equation in polar coordinates (written here in modern notation):

$$r = \frac{p}{1 + e\cos v} = \frac{a(1 - e^2)}{1 + e\cos v},$$

where $e$ is the eccentricity of the orbit or conic section, $v$ is the true anomaly angle, $a$ is the semimajor axis, and $p$ is the semi-latus rectum (half the width of the ellipse) or parameter of the orbit, giving the distance between the focus and a point of the conic section for $v = \pi/2$. The true anomaly angle $v$ defines where a body is within the orbit with respect to perigee (the point in the orbit at which a body is closest to Earth), and is the only orbital element that changes with time.

The orbital eccentricity of an astronomical object is a parameter that determines the amount by which its orbit around another body deviates from a perfect circle. The term derives its name from the parameters of conic sections, as every Kepler orbit is a conic section. In a two-body problem with inverse-square-law force, every orbit is a Kepler orbit. The eccentricity of this Kepler orbit is a non-negative number that defines its shape. For elliptical orbits, the eccentricity can be simply calculated from the periapsis (closest distance) and apoapsis (farthest distance of the orbit to the center of mass of the system, which is a focus of the ellipse) radii. In his first lunar theory,[18] Euler assumed the eccentricity of the Earth's orbit as $e = 0.01678$.[19]

In the same memoir Euler introduced an iterative method to solve Kepler's equation. Solving Kepler's equation is of great importance not only in astronomy but in astronautics as well since we need it to determine the position of artificial satellites and space probes at every point in their trajectories. Let's first review the essence of the problem to introduce Euler's work.

## Kepler's Problem

In his *Astronomia nova* and latter in *Epitome Astronomiae Copernicanae*, Johannes Kepler[20] presented an equation for the solution of planetary orbits. He introduced the eccentric anomaly $E$, and the mean anomaly $M$. The term *anomaly* means irregularity and astronomers have used it (instead of angle) to describe planetary position because the observed locations of a planet often showed small deviations from the predictions. The mathematical relationship between these anomalies and the eccentricity of the orbit is now known as Kepler's equation.

---

[17] Euler, L., *De motu planetarum et orbitarum determinatione* (On the motion of planets and orbits). E37. presented to the St. Petersburg Academy on November 21, 1735.Sect. 11, p. 72.

[18] The eccentricity of a celestial body's orbit depends on the specific orbital energy (total energy divided by the reduced mass), the standard gravitational parameter based on the total mass, and on the specific relative angular momentum (angular momentum divided by the reduced mass).

[19] Recueil pour les Astronomes par M. Jean Bernoulli, Tome III., p.327.

[20] Kepler derived the equation in 1609 in Chapter 60 of his *Astronomia nova*, and also in book V of his *Epitome of Copernican Astronomy* (1621). He proposed an iterative solution to the equation.

As shown in the figure below, the mean anomaly $M$ is the angular distance from perihelion which a (fictitious) planet, moon, or satellite would have if it moved on the circle of radius $a$ with a constant angular velocity and with the same orbital period $T$ as the actual planet moving on the ellipse. By definition, $M$ increases linearly (uniformly) with time: $M = \frac{2\pi t}{T}$ ($t$ is the time since periapse, and $T$ is the period of the motion).

Kepler's equation gives the relation between the polar coordinates of a celestial body (such as a planet or moon), or an artificial satellite or spacecraft, and the time elapsed from a given initial point. Kepler's equation is of fundamental importance in celestial mechanics, but cannot be directly inverted in terms of simple functions in order to determine where the body will be at a given time.

Derived from the second law of planetary motion, Kepler's equation is defined as follows: Let $M$ be the mean anomaly (a parameterization of time) and $E$ the eccentric anomaly (a parameterization of polar angle) of a body orbiting on an ellipse with eccentricity $e$, then

$$E(t) - e\sin E(t) = M(t)$$

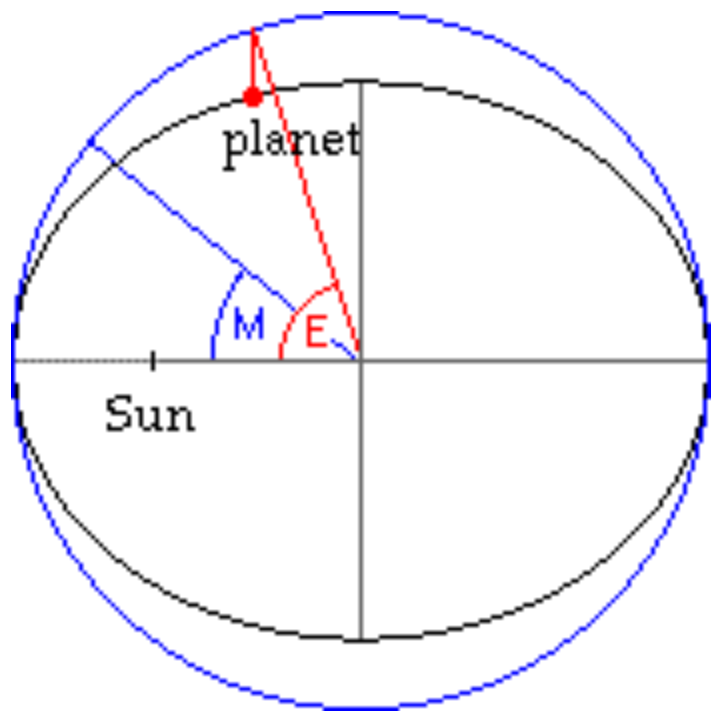

If the eccentricity $e$ and the eccentric anomaly $E$ are known, then the value of $M$ at a given time is easily found. However, usually the problem is to find $E$ when $M$ and $e$ are known. The position of the planet, moon, or spacecraft can be determined from $E$. Kepler's equation cannot be solved algebraically. We typically use numerical methods to calculate $E$, or we assume a relatively small $e$ to find an approximate solution.

In his 1735 memoir,[21] Euler considered the orbital motion of the planets in elliptical paths, with the objective of determining their positions at given times in terms of their eccentric anomaly (measured from the aphelion, the point in the orbit where a body is farthest from the Sun). He developed a recursive expression that included the eccentric anomaly (in terms of the mean anomaly) and the eccentricity. Then he continued until two successive iterations did not yield a reasonable change in the eccentric anomaly. In his calculations for the orbital position of Mars, Euler made extensive use of logarithms. He also used observational data reported by the French Academy of Sciences. However, no satisfactory solution is evident in his paper.

[21] Euler, L., *De motu planetarum et orbitarum determinatione*. E37. Presented to the St. Petersburg Academy on 21 November 1735.

In 1740, Euler amended astronomical tables[22] and, clearly not satisfied with the previous solution to Kepler's equation, he tried once again. Euler gave the solution in the form[23]

$$E = M + e(\sin M + e \sin(M + e \sin M + \cdots,$$

which reduces to the successive approximations

$$E_n = M + e \sin E_{n-1}, \qquad n = 1, 2, 3, \ldots; \quad E_0 = M.$$

For elliptical orbits ($0 < e < 1$), the iteration converges, although slowly when the eccentricity approaches one.

Nowadays, for most applications, Kepler's problem can be computed numerically by finding the root of the function, which typically is performed iteratively via Newton's method for finding successively better approximations to the roots (or zeroes) of a real-valued function.

That same year, Euler wrote a short memoir where he determined the planetary orbits around the Sun.[24] He considered the data obtained by English astronomer John Flamsteed. In 1681, Flamsteed proposed that the two comets observed in November and December of 1680 were not separate bodies, but rather a single comet travelling first towards the Sun and then away from it. Initially Newton disagreed, but eventually he had to agree with Flamsteed, theorizing that comets, like planets, moved around the Sun in large, closed elliptical orbits. Euler tried to confirm the orbit of the comet, using his own mathematical methods.

## Euler's *Mechanica* and Celestial Mechanics

Before Euler, astronomy was based on geometrical methods that were rather cumbersome. Then Euler came along and with his analytical approach, and the view of the universe took a new scientific perspective. Euler "demonstrated it by examples which, imitated since then by men of ability and reputation, may in time bestow a new form on astronomy."[25] In fact, his most important memoirs, with which Euler earned an impressive number of prizes from the Paris Academy of Sciences, relate directly or indirectly to celestial mechanics—a branch of science which demanded the highest efforts of the greatest mathematicians of the time, and which made it possible to advance astrophysics.

In the 1730s, the scientific climate in Europe brewed with ideas from two opposite thinkers: Cartesians and Newtonians, i.e., those who followed Descartes theories and those who wanted to believe Newton. With his *Principia*, Newton had stated a new principle of gravitation to explain how our Solar System works. However, Gottfried Leibniz criticized Newton for not explaining how gravity acts across space, while French scientists argued that the force of gravity was merely a supernatural idea. Many scholars had embraced René Descartes's theory of vortices, which postulated that the space is entirely filled with matter in various states, whirling about the Sun, and used the Cartesian conservation of motion, which asserts that the total quantity of motion—mass times velocity—always remains the same.

---

[22] Euler, L., *Emendatio tabularum astronomicarum per loca planetarum geocentrica* (An emmendation to astronomical tables of locating the geocenters of planets) E131. Submitted to the St. Petersburg Academy on March 28, 1740.

[23] Dutka J., "A note on Kepler's equation." Archive for History of Exact Sciences, 15 (1), pp. 59-65 (1997).

[24] Euler, L., *Orbitae solaris determinatio* (Determination of orbits around the sun). Commentarii academiae scientiarum Petropolitanae 7, 1740, pp. 86-96

[25] Eulogy to Monsieur Euler by the Marquis de Condorcet (Translated by John S.D. Glaus, The Euler Society, March 2005): http://www.math.dartmouth.edu/~euler/historica/condorcet.html

A gifted mathematician, Euler had mastered the infinitesimal calculus, which he learned from the works of Leibniz, Jacob and Johann Bernoulli. Euler had also studied the works of Galileo, Huygens and Newton, and thus he understood well the new Newtonian theories and saw how they could be improved. From the start of his career, Euler began applying the calculus for the solution of problems in mechanics.

In 1736, Euler published the first volume of his monumental book *Mechanica*.[26] In it, he restated Newton's laws and discovered new features of motion. Euler's approach was superior to Newton's, as he applied the full analytic power of calculus to dynamics. In this book, Euler recast Newton's theories of motion in a more sophisticated mathematical language, and outlined a program of studies embracing every branch of science, involving a systematic application of analysis. The two-volume book laid the foundations of analytical mechanics.[27] He was twenty-nine years old!

Using axioms, definitions, and logical deductions, Euler built a rational science of mechanics. He was the first to use partial differential equations with homogeneous functions of two independent variables (Part 2,§ 832). Euler was the first scholar to write Newton's second law in a differential form. The equation $F = ma$ first appeared in a memoir[28] that he published in 1750. Euler also introduced the moving coordinate system, consisting of tangent and normal to the trajectory, which accompanies the moving point. He used the term "infinitely small body, that is more or less a point" as a concept required by analytical mechanics, independent of all philosophical and physical speculations. Euler's contribution had a huge impact in celestial mechanics, and his works inspired many other scientists after him, especially Joseph-Louis Lagrange who stated in his *Mécanique analytique* (1788), that Euler's *Mechanica* was "The first great work where Analysis has been applied to the science of movement."

## Academy Prize for a Study of the Ocean Tides

In the spring of 1738, the French Académie Royale des Sciences offered a prize on *le flux et le reflux de la mer,* to explain the causes of the tides. This was an exciting opportunity for Euler who had until the fall of 1739 to expand on his ideas concerning Newtonian mechanics and advance the knowledge of this natural phenomenon. Tides are the rise and fall of sea levels caused by the combined effects of the gravitational forces exerted by the

---

[26] Euler, L., *Mechanica sive motus scientia analytice*, Tomus I. First volume published in 1736. (E15) Opera Omnia: Series 2, Volume 1.
[27] On Euler's contribution to classical mechanics see Truesdell's *Essays in the History of Mechanics*, Springer-Verlag, Berlin, 1968.
[28] Euler, L., *Découverte d'un nouveau principe de Mécanique* (Discovery of a new principle in Mechanics) (E177). Presented to the Berlin Academy on 3 September 1750.

Moon and the Sun and the rotation of the Earth. This phenomenon had been studied for centuries. Galileo, for example, attributed the tides to the sloshing of water caused by the Earth's movement around the Sun. Newton, on the other hand, explained the ebb and flow as the product of both lunar and solar gravitational attractions as the origin of the tide-generating forces.

Euler submitted a memoir,[29] where he showed that the horizontal component of the tidal force alone (not the vertical component as considered by Newton) was effective in raising the tide, and he discovered an exact geometrical form for the tidal potential.[30] Euler explained the tides based on the Newtonian principle of universal gravitation. He first reviewed the various explanations of tides that had been proposed, making criticisms along the way.[31] Then Euler calculated the magnitudes of the tide-generating forces and concluded the following: (a) the lunar tide is greater than the solar tide; (b) a sphere is attracted as if its entire mass is concentrated at its center.

Euler developed a theory based on geometrical and analytical considerations by assuming that the vertical motion of the water behaves as a forced harmonic oscillator with amplitude and phase influenced by the inertia of the water.[32] He then described the retardation of the tides, explaining that they are the result of the effect of the inertia of the water. Euler explained the principal phenomena of the tides in open oceans and near islands and compared them with observations.

In the spring of 1740, the French Academy of Sciences announced the winners. Euler had won the award (his second Academy prize), along with Daniel Bernoulli, Scottish mathematician Colin Maclaurin,[33] and P. Cavalleri (spelled Cavalier in the French documents).[34] Just like Euler, Bernoulli and Maclaurin had taken Newton's theory as baseline for their own dissertations.

## Vision Loss in 1738

When Euler took over the supervision of the Department of Geography, his salary was increased to twelve hundred rubles, and he began an intensive work of collaboration with Delisle on the cartography of Russia. During that time Euler suffered a nearly fatal fever and three years later, in 1738, he lost vision in his right eye —Euler was thirty-one years old! The end of Euler's 10-year astronomical observations at the St. Petersburg observatory coincides with his vision problems.

In a letter to Goldbach on 21 August 1740, Euler wrote: "Geography is fatal to me," alluding to his cartography work. Then he added: "You know that I have lost an eye and [the other] currently may be in the same danger." In the *Éloge*, his close colleague and grandson-in-law Nicolaus Fuss[35] stated: "three days of intense astronomical calculations connected with geographical work underlay the loss in the right eye and began a course leading to Euler's total blindness in 1767."

---

[29] Euler, L, *Inquisitio physica in causam fluxus ac refluxus maris* (A physical inquiry into the cause of the ebb and flow of the sea). Presented to the St. Petersburg Academy on June 15, 1739.

[30] Cartwright, David E., *The Tonkin Tides Revisited*. Notes and Records of the Royal Society of London, Vol. 57, No. 2 (May, 2003), pp. 135-142.

[31] Summary at the Euler Archive based on Eric J. Aiton's introduction (written in English) to Opera Omnia Series 2, Volume 31 and on Clifford A. Truesdell's introduction to Opera Omnia Series 2, Volume 12) . http://eulerarchive.maa.org/pages/E057.html

[32] *Ibid*.

[33] Maclaurin is known for his contributions to geometry and algebra.The Maclaurin series, a special case of the Taylor series, is named after him.

[34] P. Antoine Cavalleri was a Jesuit professor of mathematics at Toulouse. He wrote *Dissertation sur la cause physique du Flux et Reflux de la Mer*.

However, blindness did not hinder his genius and scholarly productivity. Euler was a very prolific mathematical writer and published hundreds of scientific papers. The *Opera Omnia*, published by Birkhäuser and the Euler Commission of Switzerland, contains most of Euler's works. Publication began in 1911, and to date 76 volumes have been published, comprising almost all of Euler's works. Euler wrote about 800 books and papers.[36] The "official" number of entries in Eneström's index is 866, including a number of letters and unfinished manuscripts. The Euler Archive[37] currently has 834 works. During his amazing and productive career, Euler won the Paris Academy of Sciences Prize twelve times! [38]

Euler's contribution to astronomy continued after he discontinued his daily work at the St. Petersburg observatory. For fourteen years, Euler was devoted to the general promotion of science in Russia, writing textbooks for the local schools and many important memoirs in diverse areas of pure mathematics and mathematical physics.

## 3. Invitation to Berlin by the Prussian King

While mathematicians and philosophers advanced the sciences, the political world changed and new leaders were taking power, changing the social landscape and, to some extent, affecting the lives of some prominent scholars. In Russia, for example, when Empress Anna Ioannovna died on 17 October 1740, a new regency took over and the socio-political situation in St. Petersburg became unstable.

Meanwhile, closer to Euler's homeland, the death of Frederick William of Prussia on 31 May marked the beginning of a new era in "the German lands," led by his twenty-eight year old son, Frederick II. When he took the throne of the Prussian Kingdom, one of his first actions was to revitalize the Prussian Academy of Sciences, and he invited the most prominent men of science to join it. Frederick II made French the official language of the court, including the Academy, and speculative philosophy was the most important topic of study. In the 1740s, the membership of the Prussian Academy was strong in mathematics and philosophy.

In 1740, French scholar Pierre-Louis Maupertuis went to Berlin at the invitation of Frederick II of Prussia. Maupertuis suggested to the ambitious young king to extend the invitation to Euler. Having such a distinguished young scholar would add prestige to the Berlin Academy of Sciences, as Euler was the leading mathematician of Europe, having already won two prizes from the French Academy of Sciences.

The invitation from Frederick must have come in the nick of time. Euler had felt insecure with the new Russian regime, as it imposed rules that affected his work at the

---

[35] Nicolaus (Nikolai) Fuss (1755-1826) married Euler's granddaughter, Albertine Benedikte Philippine Luise Euler (1766-1822). She was the daughter of Johann Albrecht Euler (1734-1800) and his wife Anna Sophie Charlotte Hagemeister.

[36] Sandifer, E., How Euler Did It, MAA Online (2005).

[37] The Euler Archive, http://eulerarchive.maa.org/

[38] These don't include the three awards he won with his two sons. Euler's winning essays are posted and annotated at the Euler Archive: http://eulerarchive.maa.org/publications/journals/ParisPrizePieces.html

Academy. He was aware of the political unrest and open hostilities directed at the German community. His wife Katharina was fearful of the fires that frequently broke out. In addition, Euler was experiencing health problems and personal grief—in the past four years he had lost three of his small children. Euler asked permission to leave the Russian Academy, saying: "… I find it necessary, for the sake of my weak health and for other circumstances, to seek a more pleasant climate and to accept the summons made to me by his Prussian Royal Highness. Because of this I ask the Imperial Academy of Sciences to grant me the favor of releasing me and to supply me and my household with the necessary passport for travel …"[39]

Permission was granted and in 1741, Euler became subject of Frederick the Great, King of Prussia. On June 19, Euler left St. Petersburg with his family—his wife and two children aged seven and one, and his brother Johann Heinrich—and servants. They arrived in Berlin at the end of July and immediately Euler bought a nice house, although the king had given him the privilege of a free residence near the Academy. Officially named Académie Royale des Sciences et Belles Lettres, its first session took place in January 1744 and on 12 May 1746, Frederick II appointed Maupertuis as president and Euler as Director of the Mathematical Class.

At thirty-four, Euler had still a universe of ideas that would mature in the following decades. During his first years in Berlin, Euler published his first book of astronomy and dozens of articles related to his research in celestial mechanics. Euler also established a new branch of mathematics (calculus of variations), wrote the book that lays the foundations of modern mathematical analysis, and derived new methods to tackle many problems in mathematical physics.[40]

Despite his partial blindness, Euler once again became involved with observational astronomy at the Berlin observatory. Originally, the Prussian Academy of Sciences had no such facility; then in 1711, the first small observatory was furnished for astronomer, Gottfried Kirch, financing itself through calendar computations. When Euler joined the Prussian Academy, the observatory was in poor shape following the death of his last director, astronomer Christfried Kirch. His sister, Christine Kirch, had continued the calculations for the calendars and was responsible for keeping the accounts. She also assisted astronomers in the use of the observatory.

Euler could not help but compare the Berlin observatory with that in St. Petersburg. He wrote to Schumacher praising the Russian institution: "the building in which it is housed is so well adapted to astronomical aims that we are unable to propose a better model in that respect."[41] Thus, upon arriving in Berlin, Euler spent considerable effort in reinstating and improving the Prussian astronomical facility. By the spring of 1748, Euler informed Delisle that the restored Berlin observatory was equipped with a *camera obscura* resembling the one in St. Petersburg,[42] where he carried out experiments with "artificial eclipses" using Delisle's methods.

---

[39] Shuchat A., and S. Gindikin, Tales of Mathematicians and Physicists. Springer, April 2007. p. 179

[40] A catalog of the books and papers written by Euler throughout his career available at The Euler Archive.

[41] Nevskaya, N.I. and K.V. Kholshevnikov, Euler and the Evolution of Celestial Mechanics, in Euler and Modern Science, eds. N.N. Bogolyubov, G.K. Mikhailov, and A.P. SYushkevich . MAA Vol. IV. p. 284.

[42] *Ibid*. p. 285.

His fascination with celestial mechanics is evident. Shortly after arriving in Berlin, Euler presented to the Academy the results of his work on the motion of comets,[43] showcasing the analysis performed to determine the orbit of the comet observed by Cassini in 1742. In 1744, he presented his New Astronomical Tables for Calculating the Position of the Sun.[44] These tables were of outmost importance to celestial navigation. This is a position fixing technique that has evolved over several thousand years to help sailors cross the oceans without having to rely on estimated calculations, or dead reckoning, to know their position. Celestial navigation uses sights, or angular measurements taken between a celestial body and the visible horizon. The Sun is most commonly used, but navigators can also use the Moon, a planet or one of 57 navigational stars whose coordinates are tabulated.

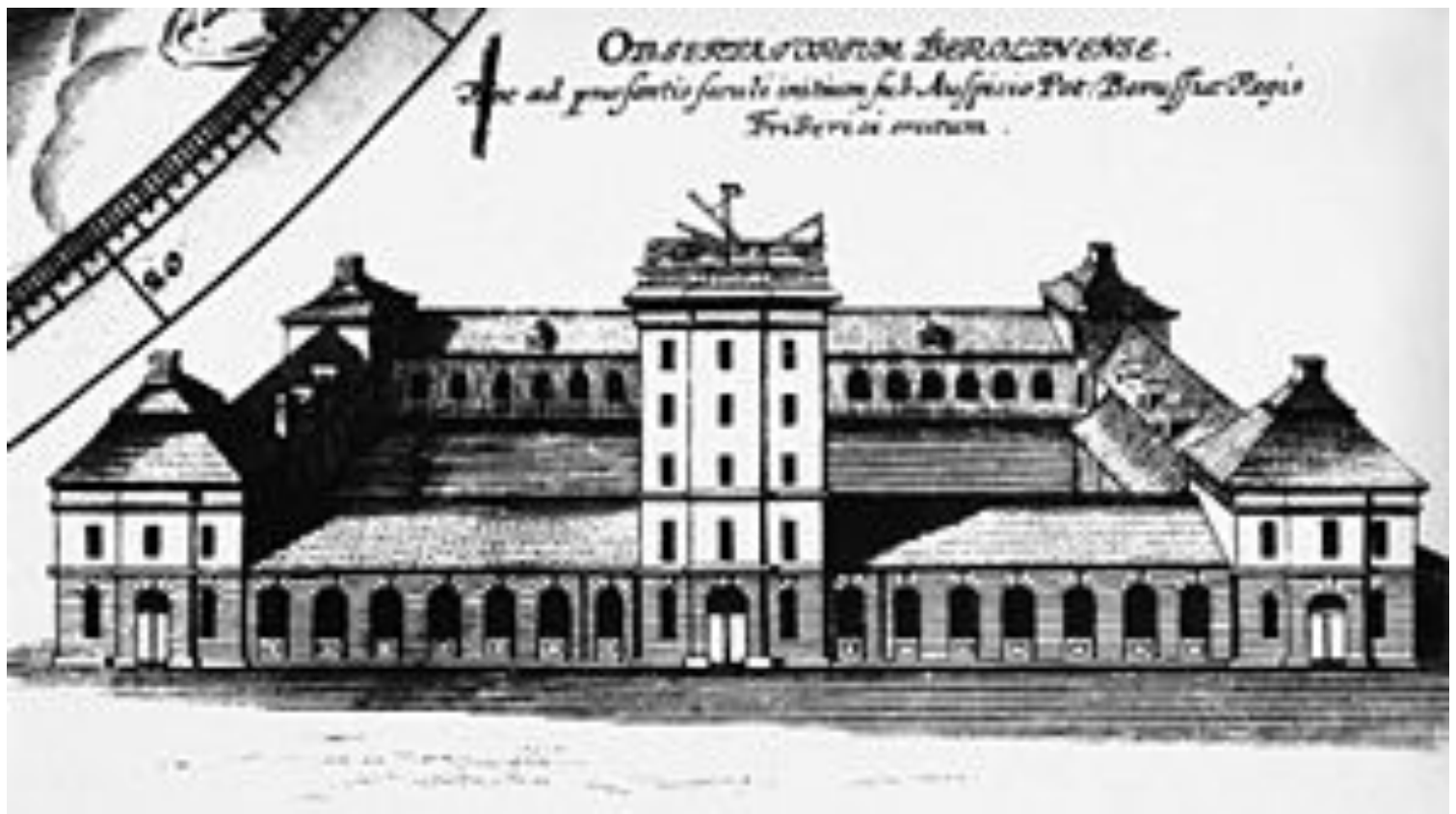

First Berlin Observatory in Dorotheenstraße (Berlin, 1711). Source: Detail der Tafel 19, die südliche Himmelssphäre, aus dem Atlas Coelestis von Johann Gabriel Doppelmayr, (1742, Homännische Erben, Nürnberg), rechts unten: Die alte Berliner Sternwarte. Vermutlich Scan eines Exemplars im Besitz des Leibniz-Instituts für Astrophysik Potsdam

### First Book of Astronomy

Shortly after arriving in Berlin, Euler began to write "Theory of the Motions of Planets and Comets,"[45] his first full-length astronomy book. In this highly didactic text, he introduced new methods of investigating planetary perturbations. In a clear and simple manner, Euler stated specific problems, providing their solutions aided by scholia and corollaries. He recognized that the orbit could be shaped as a circle ($e = 0,\ p = r$),

---

[43] Euler, L., *Determinatio orbitae cometae qui mense Martio huius anni 1742 potissimum fuit observatus* (Determination of the motion of a comet which can be observed in March of this year, 1742). E58. Presented to the Berlin Academy on September 6, 1742.

[44] Euler, L., *Nouvelles tables astronomiques pour calculer la place du soliel* (New astronomical tables for calculating the position of the sun). Read to the Berlin Academy on April 9, 1744, but only the abstract and the tables appeared in the Hist. de l'acad. d. sc. de Berlin, (1745), 1746, pp. 36-40 + 1 table; the presentation date is found on p. 36 of the same [E836a].

[45] Euler, L., *Theoria motuum planetarum et cometarum* (Theory of the motions of planets and comets). Published in Berlin in 1744

ellipse ($0 < e < 1$; $p = a(1 - e^2)$), parabola ($e = 1$; $p = p$) or hyperbola ($e > 1$; $p = a(e^2 - 1)$. The names of the curves derive from the fact that they can be obtained as the intersections of a plane with a right circular cone. If the plane cuts across a half-cone, the section is an ellipse; one obtains a circle, if the plane is normal to the axis of the cone, and a parabola in the limit case of a plane parallel to a generatrix of the cone.

Euler expressed the polar equation of the orbit as

$$r = \frac{ab}{a + (b - a)\cos v}$$

where $a$ is the perihelion distance. We recognize this form as the conic section describing a parabola. Parabolic orbits are open orbits but contain a special boundary case with $e = 1$. These orbits describe the path of a body falling in towards another from very far away when it starts with very low velocity. Long period comets follow this type of orbit.

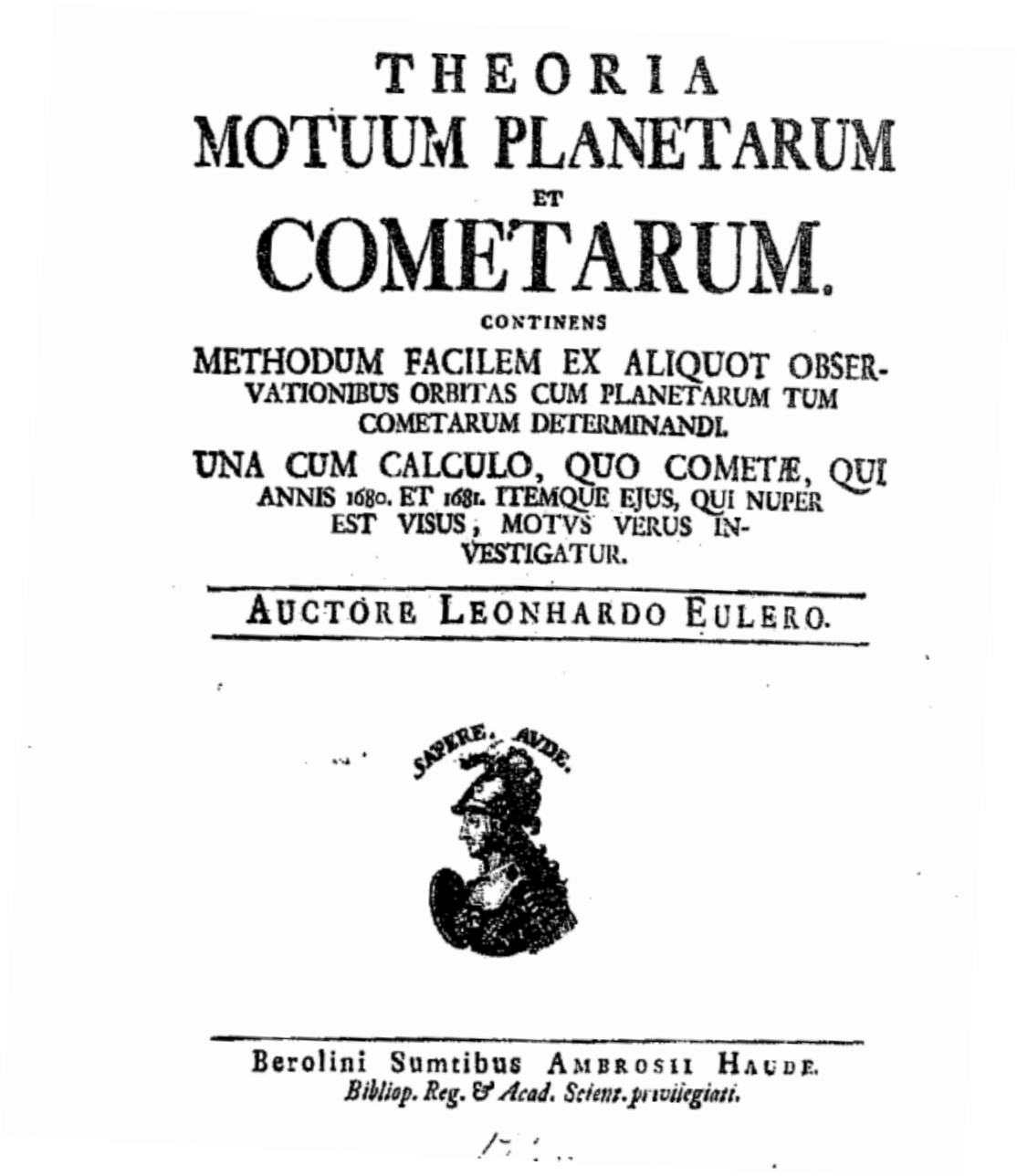

THEORIA
MOTUUM PLANETARUM
ET
COMETARUM.
CONTINENS
METHODUM FACILEM EX ALIQUOT OBSERVATIONIBUS ORBITAS CUM PLANETARUM TUM COMETARUM DETERMINANDI.
UNA CUM CALCULO, QUO COMETÆ, QUI ANNIS 1680. ET 1681. ITEMQUE EJUS, QUI NUPER EST VISUS, MOTVS VERUS INVESTIGATUR.
AUCTORE LEONHARDO EULERO.
SAPERE AUDE.
Berolini Sumtibus AMBROSII HAUDE.
*Bibliop. Reg. & Acad. Scient. privilegiati.*

Euler's Book "Theory of the Motions of Planets and Comets" (1744).

For his theory of motion of celestial bodies, Euler derived differential equations based on the principles of mechanics.[46] In this paper, Euler discussed the observed irregularities of the planets as they orbit the Sun. He obtained the solution, which is said to be a regularization of the inverse problem of Newton. His analysis was based on the laws of linear and angular momentum. In fact, Euler was the first to consider Newton's laws as general principles that are applicable to each part of every macroscopic system. He introduced the idea that forces are three-dimensional vectors, the notion of reference frame, and he used rectangular Cartesian coordinates to define the motion of the celestial bodies (planets in the Solar System). Euler introduced for the first time the Newtonian equations in the form that we recognize today. For example, he defined the velocity of a body as given by its three components $\frac{dx}{dt}, \frac{dy}{dt}, \frac{dz}{dt}$, and the acceleration as $\frac{d^2x}{dt^2}, \frac{d^2y}{dt^2}, \frac{d^2z}{dt^2}$.

[46] Euler, L., *Recherches sur le mouvement des corps célestes en gènèral* (Studies on the movement of celestial bodies in general, E112. Presented to the Berlin Academy on 8 June 1747.

In 1748, the Berlin Academy published the French translation of Euler's memoir explaining Kepler's law of planetary motion.[47] He began by establishing how the planets move along the ecliptic and explained retrograde motion, that is, when the planets exhibit a very irregular motion when observed from the Earth. "However, if we observe the planet from the stand point of an observer on the Sun, this retrograde motion will not occur, and only a west to east path of the planet is seen."[48] At the end, Euler provided a table that would be useful to astronomers looking for the eccentricity of a planet's orbit.

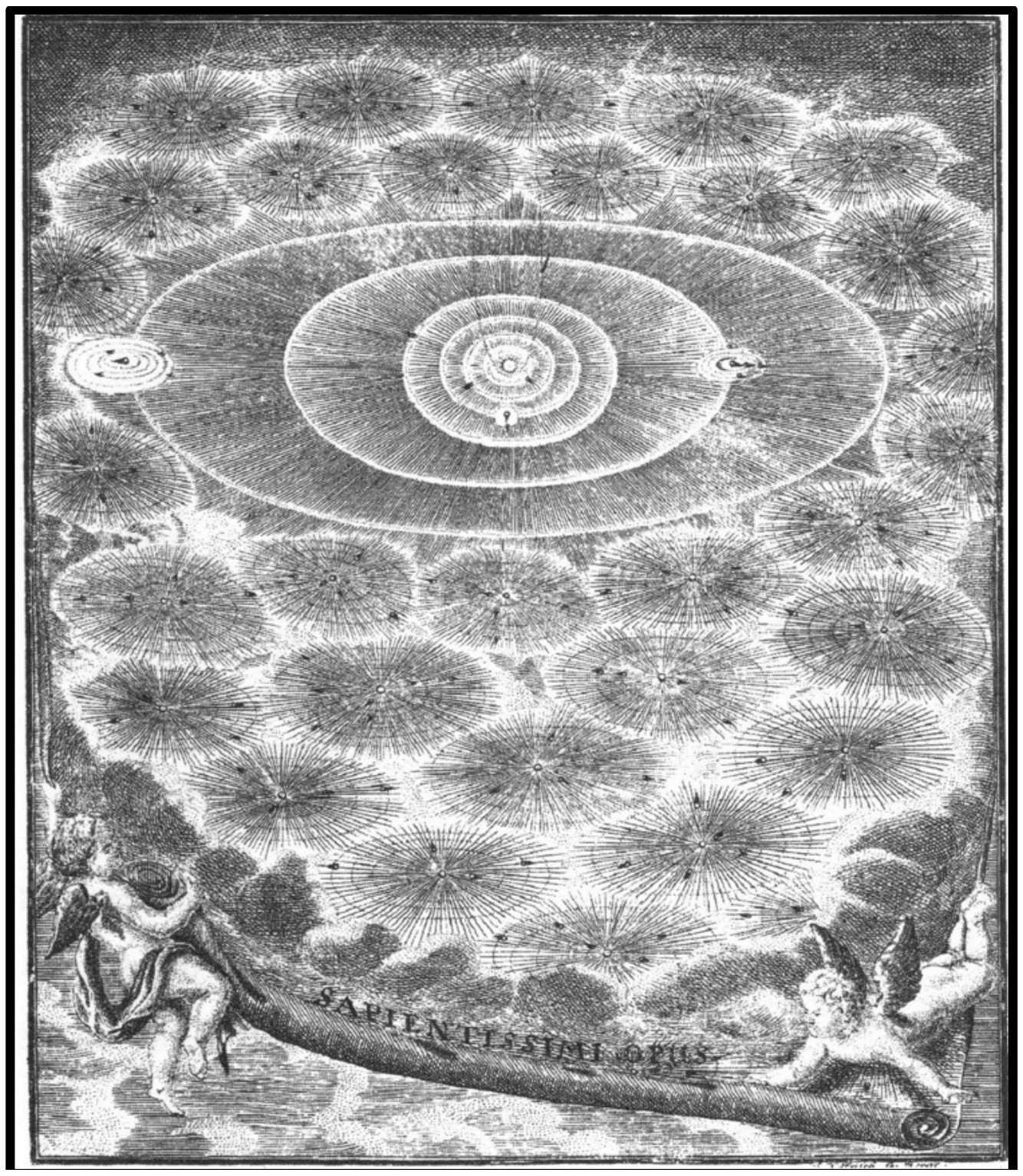

Beautiful illustration in Euler's 1744 "Theory of the planets and comets." Engraving depicts the solar system with the then known planets, including Jupiter and Saturn with their four (known at the time) moons. The trajectory of the comet is indicated with dotted line passing between the Sun and Mercury. Frontispice allégorique illustrant le mouvement des planètes. [Cote :2198A]

On 25 July 1748, a solar eclipse was visible in Berlin, and despite his vision problems, Euler made arrangements to study it, assisted perhaps by Christine Kirch.[49] His former colleague Delisle, now in Paris, had sent numerous Avertissements (information bulletins) to astronomers, encouraging them to observe the solar eclipse and record their impressions.

---

[47] Euler, L., *Memoire sur la plus grande equation des planetes* (Memoir on the Grand Equation of the Planets). Memoires de l'Académie des sciences de Berlin 2, 1748, pp. 225-248
[48] Synopsis of Euler's paper by Thomas J. Osler and Jasen Andrew Scaramazza of Rowan University.
[49] Sandifer, E., How Euler Did It: *Fallible Euler*, MAA Online, February 2008.

In his 1748 essay,[50] Euler reported: “En observant les moments de l’Eclipse du Soleil que nous eûmes ici le 25 Juillet 1748. & en rendant compte de mes observations, je n’avois eu pour but que d’arriver à une détermination plus exacte du véritable mouvement de la Lune & de sa parallaxe.” During this astronomical event, Euler was not just a casual observer but clearly was working as an astronomer taking his own data. He added: “Mais il n’a pas laissé de se présenter dans le cours de l’Observation de cette Eclipse quelques autres Phénomènes remarquables, qui ne dépendaient, ni du mouvement de la Lune, ni de la parallaxe, mais qui semblaient donner à connaitre la réalité de la réfraction des rayons qui rasent les bords de la Lune, & décider cette question agitée depuis long temps parmi les Astronomes; Si la Lune est environnée d’une Atmosphère, ou non? C’est ce qui m’engage à examiner ici plus exactement les phénomènes de cette espèce, que j’ai observés pendant cette Eclipse, & à en recherches les causes ...” Euler’s objective was to answer the question asked by many astronomers: is the Moon surrounded by an atmosphere, or not?

Euler wanted to demonstrate that certain phenomena that resulted from the eclipse were evidence that the Moon had an atmosphere and that it was almost 200 times less dense than that of the Earth.[51] He was wrong, of course. More than two hundred years later, the American astronauts who landed on the Moon proved that it has an atmosphere so tenuous as to be nearly a vacuum. In fact, David R. Scott, the seventh astronaut to walk on the lunar surface, simultaneously dropped a hammer from one hand and a feather from the other. In the absence of atmospheric friction both hammer and feather fell down at the same time.

## Lunar Theory: The Restricted Three-Body Problem

In the early eighteen century, Newton’s gravitational theory was in question because it appeared to be inconsistent with the astronomical observations of the Moon and the motion of the planets in the Solar System. The motion of the Moon in particular, calculated by the leading mathematicians, yielded results that did not match the astronomical observations. The discrepancy was too huge to ignore.

Many anomalies or irregularities in the Moon’s motion around the Earth were known since the time of Ptolemy. Newton attempted an explanation of the cause for those irregularities, and in 1702, he published a lunar theory, which allegedly was based on his gravitational theory, although the analysis used epicycles (a Ptolemaic geometric model used to explain the variations in speed and direction of the apparent motion of celestial bodies). Newton issued rules for finding the position of the Moon, but he did not justify the rules with gravitational theory.[52] Newton’s lunar theory was far from adequate.

---

[50] Euler, L., *Sur l’atmosphere de la Lune prouvée par la derniere eclipse annulaire du Soleil*. (On the atmosphere of the moon as proved by the last ringed eclipse of the sun), E142. According to C. G. J. Jacobi, a Latin treatise with the title: “*De atmosphaera lunae ex eclipsi solis evicta*” was presented to the Berlin Academy on December 5, 1748.

[51] The Euler Archive, http://eulerarchive.maa.org/

[52] Kollerstrom, N., *Newton’s 1702 Lunar Theory*, Department of Science & Technology Studies, University College London.

Moreover, his analysis was not sufficiently accurate to predict the mass of the Moon. In the first edition of his *Principia*, Newton overestimated the lunar mass by a factor of three, finding that: "... the mass of the Moon will be to the mass of the Earth as 1 to 26, approximately," citing the relative densities as 9 to 5. The actual ratio is 1:81, as the lunar density is merely 0.6 that of Earth.[53]

In 1747, after attempting to resolve the problem of the lunar orbit, French mathematician Alexis Clairaut concluded that Newton's theory was incorrect. And since d'Alembert's calculations agreed with his, Clairaut proposed a correction to Newton's inverse square law of gravitational attraction by adding a term. However, a few months later Clairaut realized that the difference between the observed motion of the Moon's apogee and the one predicted by the theory was due to mathematical errors in the approximations rather than Newton's law.

Attempting to take the lunar theory beyond that stage left by Newton, in 1746 Euler published his first lunar tables.[54] These tables, based on his theoretical analyses, were designed to facilitate the calculation of planetary positions, lunar phases, eclipses and provide other data for calendars. However, due to inaccurate coefficients (derived from other astronomers' observations), Euler's tables were still inadequate. In the years that followed, Clairaut and D'Alembert computed more accurate tables. By then Euler had formalized the restricted three-body problem, while attempting to solve analytically—for the first time—the general problem of orbital perturbations.

Euler investigated the slowing down of planetary motion, which he attributed to the resistance of the (hypothetical) aether. At that time he, like most scientists and philosophers, believed that the aether was a kind of subtle fluid that filled the spaces between material substances and provided a medium for transmitting phenomena like gravity, light, magnetism and electric fields. Euler hypothesized that the aether might exert friction affecting the motion of bodies. In 1746, he wrote one article[55] and letters to Delisle on the resistance the aether exerted on the motion of planets and comets in their orbits around the Sun.[56]

It was natural for Euler to accept such theory. The notion of the aether may have been as intriguing as the idea of dark matter/dark energy that contemporary physicists use to explain the expansion of the universe. Years ago many scientists believed that the expansion of the Universe was slowing due to gravity. However, observations of distant galaxies by the Hubble Space Telescope in 1998 indicated that the expansion has been accelerating! No one expected this, and no one knew how to explain it. Some thought that the discrepancy was a result of a long-discarded version of Einstein's theory of gravity that contained what was called a cosmological constant. Others thought that there may be some strange kind of energy-fluid that fills space; something we cannot see and have yet to detect directly is giving galaxies extra mass, generating the extra gravity they need to stay intact. This strange and unknown matter was called "dark matter" since it is not visible.[57] To date, scientists still don't know what is the correct explanation for the

---

[53] *Ibid*.
[54] Euler, L., *Tabulae astronomicae solis et lunae* (Astronomical table of the sun and the moon). E
[55] Euler, L., *De relaxione motus planetarum* (On the enlarging of the motion of the planets).
[56] Sandifer, E., How Euler Did It: Fallible Euler, MAA Online, Feb 2008.
[57] For a brief description of the nature of dark matter/dark energy, see http://science.nasa.gov/astrophysics/focus-areas/what-is-dark-energy/

cosmic acceleration, but many papers and books are being written explaining it in terms of dark matter and dark energy. In the case of the aether, a hundred years after Euler the experiments performed by Michelson and Morley in 1887 refuted the existence of such matter. Although some contemporary scientists now relate the aether to modern field theories

Nonetheless, it must have been clear to Euler that the Moon's orbit around the Earth is elliptical, with a substantial eccentricity $e$ = 0.0549 compared with Earth's.[58] Moreover, the tidal effect of the Sun's gravitational field increases the eccentricity when the orbit's major axis is aligned with the Sun-Earth vector or, in other words, when the Moon is full or new. Euler presented a "*Theoria motus lunae*" to the Berlin Academy on 22 April 1751. The Lunar Theory is an explanation by mathematical reasoning of perturbations in the movements of the Moon founded on the law of gravitation.

Unsure about Clairaut's results, Euler had the St. Petersburg Academy establish a prize in 1752 for a theory of the Moon. Clairaut submitted a memoir describing his mathematical analysis (*Théorie de la lune*), which Euler praised. However, Euler must have been frustrated for not obtaining himself the required solution to fully account for the observed motion of the Moon.

A year later, Euler published his own lunar theory in a voluminous 350-page long book[59] divided into eighteen chapters, plus an addendum where he included his research on the inequalities of the lunar motion. The foreword (written by an unknown person) indicates that the winning essay of Clairaut, about which Euler sent a report to the Russian Academy, had inspired this work. Euler had developed a fundamental method for the approximate solution of the three-body problem.

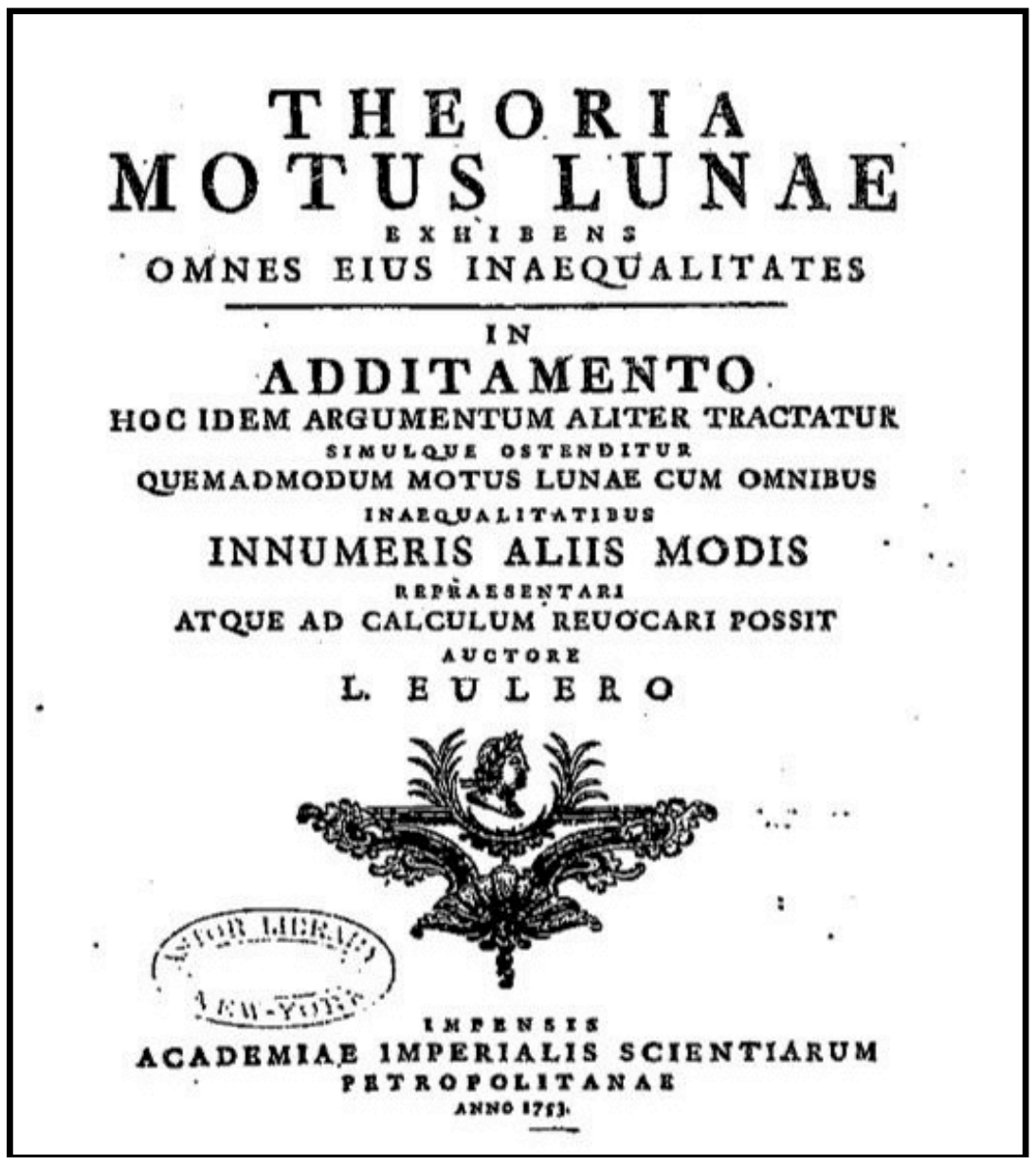

THEORIA
MOTUS LUNAE
EXHIBENS
OMNES EIUS INAEQUALITATES
IN
ADDITAMENTO
HOC IDEM ARGUMENTUM ALITER TRACTATUR
SIMULQUE OSTENDITUR
QUEMADMODUM MOTUS LUNAE CUM OMNIBUS
INAEQUALITATIBUS
INNUMERIS ALIIS MODIS
REPRAESENTARI
ATQUE AD CALCULUM REUOCARI POSSIT
AUCTORE
L. EULERO

IMPENSIS
ACADEMIAE IMPERIALIS SCIENTIARUM
PETROPOLITANAE
ANNO 1753.

Frontispiece of Euler's First Lunar Theory (1753).

[58] Earth's orbital eccentricity is $e$ = 0.017. The mean Moon-Earth distance is 384,401 km.
[59] Euler, L., *Theoria motus lunae exhibens omnes eius inaequalitates* (Theory of the motion of the moon which exhibits all its irregularities).E187. Published as a book in 1753. Opera Omnia: Series 2, Volume 23, pp. 64 - 336

The Moon-Earth-Sun system is a complicated three-body problem. The Moon's orbital inclination, combined with the inclination of the Earth's axis of rotation, causes the Moon's declination (as observed from the Earth) to vary between ±28.5° when the Moon's inclination adds to that of the Earth, and ±18° when the two inclinations oppose one another. The maxima and minima of declination repeat every 18.6 years, the period in which the ascending node of the lunar orbit precesses through a full circle.

In the sixteenth and seventeenth centuries, the prediction of longitude at sea was the most pressing scientific problem, and the finding of lunar longitude seemed to be a promising approach. Thus, Euler's lunar theory had an important practical consequence. In 1714 the British parliament had offered a considerable monetary prize for the determination of the longitude at high sea to within an error margin of one half degree. Decades later, the German astronomer Tobias Mayer used Euler's formulae to produce lunar tables, which made it possible to determine not only the position of the Moon, but also allowed the determination of the geographic longitude of a ship at high sea, as accurate as required by the British prize.

In 1765, a monetary award was given to Mayer and Euler: Mayer's widow received 3,000 pounds, and Euler 300 pounds for the theory underlying Mayer's tables. These lunar tables became part of all navigational almanacs and served merchant shipping for more than a century. Eric J. Forbes translated the Euler-Mayer correspondence and provided English readers for the first time, in 1971, compelling material to examine the influence Euler had on Mayer's work.

## Theory of Jupiter and Saturn

At the time of Euler, the two largest planets known were Jupiter and Saturn. Questions about the inequalities in their motions had been unsolved for centuries. In 1625, Kepler wrote that the locations of the two planets were more or less advanced as they should be when their mean motions were determined according to ancient observations of Ptolemy and compared with those of Tycho.

After Newton, astronomers noted that their observations showed inequalities in the motions of the giant planets that could not be accounted for by the theory of gravitation, since it suggested that the mutual action should be detectable around the time of conjunction. By the 1740s, researchers were certain that the periodic inequalities in the motion of Saturn depended on its position relative to Jupiter. In 1746, French astronomer Pierre-Charles Le Monnier wrote that not only there are some periodic inequalities on Saturn that depends on its position relative to Jupiter, "but in the same configurations returning after fifty-nine years, the error in the Tables is always growing."

At the urging of Le Monnier, in 1748 the Paris Academy established a prize competition to develop "a theory of Jupiter and Saturn explicating the inequalities that these planets appear to cause in each other's motions, especially about the time of their conjunction."

Euler was the first to formulate a theory of Jupiter and Saturn, aiming to match the observations. For his work Euler won the 1748 Paris Academy prize. However, Euler's application of his solution to observations spanning the years from 1582 to 1745 left some errors.[60] The Academy reissued the same topic for the prize contest of 1750 and several mathematicians, including Euler and Daniel Bernoulli, submitted their theories, but no prize was awarded and the problem remained unsolved. Various hypotheses were introduced, such as a resisting medium in space, and a finite time for the action of the attracting forces, but none of these were satisfactory.

Finally, for the 1752 contest, Euler developed a beautiful method of the variation of constants, introducing trigonometric series to address the question of the inequalities in the movements of Saturn and Jupiter. In the first pages, Euler discussed this calculus of trigonometric functions and showed their importance in analysis. Despite the strength of the mathematical method, the resulting variations in the mean motion of Jupiter and Saturn were contrary to the astronomical observations. Euler's results were wrong. However, by studying the perturbations on the large planets, Euler laid down the foundations of the mathematical method of perturbations, which was later refined by Lagrange. In 1786, Laplace discovered that the inequality between Jupiter and Saturn arose from the small divisor introduced by integration, one that depends on the relation that five times the mean motion of Saturn is nearly equal to twice the mean motion of Jupiter. Although the numerator of this coefficient is of the third order of the eccentricities, this inequality amounts to 20′ in the longitude of Jupiter, and to 48′ in that of Saturn.[61]

## Precession of the Equinoxes

Euler considered the problem of precession and nutation (a small irregularity in the precession of the equinoxes) and obtained a correct solution. The precession of the Earth's equinoxes is a phenomenon that has been known for centuries. Precession refers to a change in the direction of the axis of a rotating object. In the case of the Earth, precession means that its axis is not stationary but instead traces out a large circle with respect to the fixed stars. Discovery of the precession of the equinoxes is generally attributed to the ancient Greek astronomer Hipparchus (ca. 150 BC), though the difference between the sidereal and tropical years was known to Aristarchus of Samos much earlier (ca. 280 BC). Isaac Newton suggested a cause for the precession of the equinoxes, but he was unable to derive a mathematical theory that would lead to a quantifiable solution. Euler did.

In his 1750 essay,[62] Euler considered nine different problems to study the precession of the equinoxes, starting with defining the rotational motion of the Earth caused by the forces of the Sun and the Moon. He calculated the annual precession of the pole and also of the equinoxes as $24\frac{1}{3}$ seconds. For problem 8, Euler considered an arbitrary given time to determine the longitude and latitude of the north pole of the Earth, measuring the longitude from a given fixed star. He remarked that "Having properly established the true quantity of the mean precession of the equinoxes by comparing ancient and modern observations, it will be easy to determine, for any given time, the mean longitude of all the fixed stars, given that we have properly determined the longitude once. I will suppose therefore that we know the mean longitude of a star at a given time and that we want to determine the true longitude of it for the same time. This will be done by the method of two equations that the following two tables will supply."[63] In the last problem, Euler determined the true quantity of the precession of the equinoxes during the space of a given year.

---

[60] *Planetary Astronomy from the Renaissance to the Rise of Astrophysics*, Eds. René Taton and Curtis Wilson, Cambridge University Press (August 25, 1995).

[61] Hill's Theory of Jupiter and Saturn, The Annals of Mathematics, Vol. 5, No. 6 (1891), pp. 229-230.

[62] Euler, L., *Recherches sur la precession des equinoxes st sur la nutation de l'axe de la terre* (Research concerning the precession of the equinoxes and of the nutation of the Earth's axis). E171. Presented to the Berlin Academy on March 5, 1750.

Euler's memoir provoked the fury of French mathematician Jean le Rond D'Alembert. He argued that he had previously discovered some of Euler's results for which he had not received due credit. It is true that d'Alembert had published in June of 1749 the first correct derivation of the precession of the equinoxes in his book-long essay *Recherches sur la Precession des Equinoxes*. It is possible that Euler did not read d'Alembert's work because, about the same time that it appeared in Paris, he and his wife were dealing with tragedy—in August they lost twin babies. However, being a noble man, Euler credited d'Alembert for the work he did by publishing an *Avertissement au sujet des recherches sur la precession des equinoxes* (Notice on the subject of Research on the precession of the equinoxes), which appeared in *Memoires de l'academie des sciences de Berlin* 6, 1752.

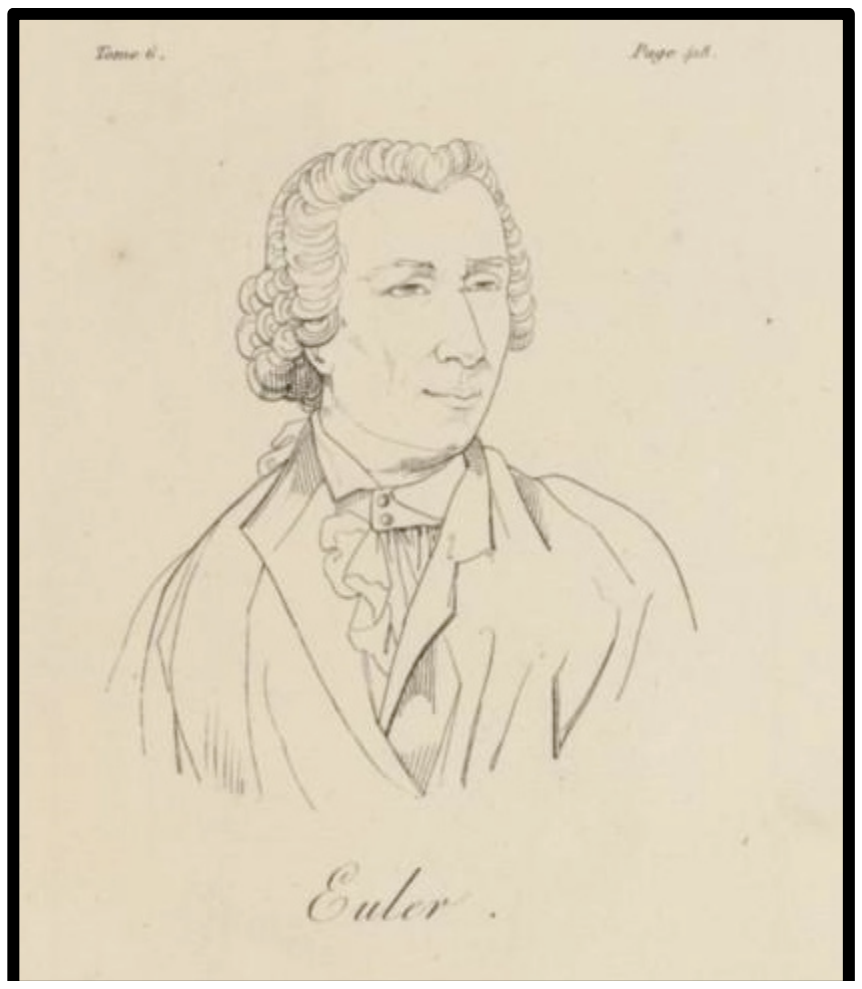

Undated Portrait of Euler. *Source*: gallica.bnf.fr / Bibliothèque nationale de France

In Berlin, Euler was fairly wealthy, owning a home near the Academy and a country estate in Charlottenburg with a great amount of corn and meadowland.[64] He and his wife had seven more children, four of which died in infancy in the years 1749-1750. In June 1750, Euler made a trip to Frankfurt a.M. to get his widowed mother so she could live in comfort with his family in Berlin. By 1765, his household has grown to sixteen, including his grown sons Johann Albrecht (31), a mathematician, Karl Johann (25), and Christoph (22), plus his brother and his family. His visual faculty continued to decline due to cataracts in his good (left) eye.

[63] Translation of Euler's article by Steven Jones and Rob Bradley.
[64] Fellmann, E.A., Leonhard Euler, Birkhäuser Verlag, Switzerland (2007).p. 101

It has been written that Frederick II was unkind to Euler.[65] The king had entrusted Euler with many duties—he had to oversee the Academy's observatory and botanical gardens, deal with personnel matters, attend to financial affairs such as the sale of almanacs, which constituted the major source of income for the Academy, not to speak of a variety of technological and engineering projects. Moreover, even though Euler was practically running the Berlin Academy, when Maupertuis died in 1759 the king refused him the title of President. The king was unable to appreciate the magnitude of the scientific contributions that came out of Euler's work. In his twenty-five years in Berlin, Euler published over three hundred memoirs (some were prize-winning essays) and books, and established new branches of mathematics. Thus, it is not surprising that Euler would feel unappreciated and decided to leave Prussia to answer the call from St. Petersburg.

Before he left, Euler recommended Italian mathematician Joseph-Louis Lagrange for the director of mathematics position that he vacated. The thirty-year old Lagrange had been communicating with Euler since 1754, sharing his brilliant mathematical discoveries. Euler was impressed with the young mathematician's method of maxima and minima, an area that both independently pioneered. The feelings were mutual, and through the years Euler and Lagrange maintained a closed correspondence, challenging one another to refine their mathematical ideas and, in several occasions, they both competed for the same honors.

Euler left in Berlin another young admirer and supporter, Johann III Bernoulli. He was the prodigy grandson of Johann I (Euler's mentor) who in 1764 had been invited by Frederick II to reorganize the observatory. Johann III was Astronome Royal in the Prussian court, and in 1771, he published the *Recueil pour les Astronomes*, a 3-volume collection of astronomy results that featured the most important work of Euler, Flamstead, Mayer, and other European astronomers. This publication also included astronomy tables, and selected memoirs of the Berlin Academy.

## 4. Return to Russia in 1766

In Russia, the winds of change were blowing in a favorable direction. When Catherine II ascended the throne in 1762, she invited Euler to return to St. Petersburg. Through the years, he had maintained close ties with the Russian Academy, and thus such invitation may have felt as a call to return home. In 1766, Euler moved to St. Petersburg, where he was given a warm welcome. He was fifty-nine, but his mind was vigorous, full of ideas that were yet to give fruit.

In addition to many other mathematical works, Euler continued solving problems of astronomy. On 8 August 1769, French astronomer Charles Messier discovered a bright comet. Messier's observations were used by astronomers all over the world to derive the comet's orbital elements. Euler was among those astronomers. In the words of Messier: "Many calculations have been made of this Great Comet by Leonard Euler and Lexell, for finding its true elliptical orbit and its period."[66] In the *Recueil pour les Astronomes*[67] published in 1771 by Johann III Bernoulli, it reports the orbital elements of the comet of 1769 calculated by Lalande, Euler, Lambert, and three other astronomers. In fact, from three selected observations, Euler and Lexell found the revolution of this comet could be 449 years to 519 in its extremes, supposing 1 [arc] minute of error in these 3 observations.

---

[65] *Ibid*. p.

[66] Messier, C., *Notice de mes comètes* (Notes on my comets). This document, written by Messier in old age by his own hand, describes his observations of all 44 comets he observed between 1758 and 1805. The Notes also contains valuable historical and more general information on Messier's time.

Euler had an extraordinary memory and the ability to perform complex mathematical analysis and computations in his head without the benefit of pencil and paper. His enormous power of memorization was paired with a rare power of concentration—noise and hustle in his immediate vicinity rarely disturbed him in his mental work.[68] These traits proved to be crucial for accomplishing his brilliant scholarly work later in life after he lost his sight in both eyes.

He developed a new method of reducing the motions of the planets in the form of astronomical tables. Already blind, Euler worked out his second lunar theory and published it in book form as *Theoria motuum lunae*,[69] Astronomer George William Hill used Euler's theory to develop his theory of Jupiter and Saturn.[70] Euler continued refining his lunar theory with the objective of carrying it to a higher level of perfection.[71] In 1862, Euler's *Opera Postuma* contained three more essays where he amended the lunar tables using observations of a lunar eclipse.

By 1771, after complications from an operation to remove a cataract in his good eye, Euler was completely blind. Nevertheless, his mathematical ideas were flowing faster than he could write them down on a large slate, and so he sought a trained mathematician to help him. In 1773, he hired Niklaus Fuss from Basel. Undaunted by blindness and other personal hardships (he lost his house to fire that year), Euler continued to publish his results by dictating them to Fuss and other assistants. Interestingly, Euler wrote the largest number of memoirs in his sixties, when he experienced some of the worst personal and health problems.

In the summer of 1778, Johann III Bernoulli traveled from Berlin to Saint Petersburg. During his visit, Johann wrote about Euler and how he coped with blindness: "... it is true the he cannot recognize people by their faces, nor read black on white, nor write with pen on paper; yet with chalk he writes his mathematical calculations on a blackboard very clearly and in rather normal size; these are immediately copied by one of his adjuncts, Mister Fuss and Golovin (most often the former) into a large book, and from these materials are later composed memoirs under his direction."[72] In this manner, Euler wrote as many articles as when he had full eyesight.

---

[67] Recueil pour les Astronomes par M. Jean Bernoulli, Astronome Royal, Tome I., p. 225.

[68] Fellmann, E.A., Leonhard Euler, Birkhäuser Verlag, Switzerland (2007).p. xx

[69] Euler, L., *Theoria motuum lunae* (Theory of lunar motion) . E418. presented to the St. Petersburg Academy on October 20, 1768.

[70] Fellmann, E.A., *Leonhard Euler*, Birkhäuser Verlag, Switzerland (2007).

[71] Euler, L., *De Theoria Lunae ad maiorem perfectionis gradum evehenda* (On the Theory of the Moon which is to be carried to a higher level of perfection) E504. Presented to St. Petersburg Academy on 6 March 1775.

[72] Fellmann, E.A., Leonhard Euler, Birkhäuser Verlag, Switzerland (2007) p. 119.

## 5. Euler and the Three-Body Problem

The three-body problem involves the mathematical description of the motion of three masses (thought of as points) subject to gravitational forces. Newton had proved that the two-body problem was integrable, and so could be solved exactly. Two-body systems can be, for example, the Earth and the Moon, or the Sun and the Earth. With Newton's formula, we tell precisely where the Earth is in its orbit around the Sun for all time, no matter how far in the past or future—as long as we assume that the Earth and the Sun are the only bodies in the system.

However, when we consider three bodies, say for example Moon-Earth-Sun, the problem is much more difficult, as the changing gravitational pull of all three moving bodies produces a very complex dynamical system. Euler realized that the three-body problem is not integrable, which means that no exact solution can ever be found. In his 1747 memoir addressing his study of motion of celestial bodies,[73] Euler formalized the restricted three-body problem (*problème restreint*), usually ascribed to Carl Gustav Jacob Jacobi and Henry Poincaré. The restricted three-body problem assumes that one of the three masses is negligible, and so exerts no gravitational influence on the other two.

When he derived his first lunar theory[74] in 1751, Euler developed a fundamental method for the approximate solution of the three-body problem. Even making certain simplifications, he realized how difficult it was to solve this problem analytically (a general close solution is impossible). First, Euler considered special cases, solving them by means of a method now called "regularization," showing that his regularized (derivatives) solution is stable with respect to variations of the data of the problem.

Euler was the first to undertake analytically the general problem of perturbation.[75] As we said earlier, he won the Paris Academy Prize for his work on the three-body problem, focusing on the Jupiter-Saturn-Sun system.[76] It was also during this time that the term "three-body problem" (*Problème des Trois Corps*) began to be commonly used by Euler and his contemporaries.

In 1762, Euler reconsidered the Earth-Moon-Sun system.[77] In this work he identified two special regions in space in the orbital configuration where a smaller mass affected only by gravity can theoretically be part of a constant-shape pattern with two larger masses (these regions are now identified as L1 and L2, shown in the figure below). The L1 point lies on the line defined by the two large masses $m_1$ and $m_2$, and between them. That is where the gravitational attraction of $m_2$ partially cancels the attraction of $m_1$. Drawing a line connecting the Earth and the Moon, there are three points where a small body may exist in an unstable equilibrium if it has a velocity giving it the same period as the Moon. In these points, the sum of the gravitational forces from the Earth and the Moon supplies just the centripetal force necessary for a circular orbit with the period of the Moon. If we visualize the gravitational forces exerted by the Earth and the Moon on a spacecraft at L1, for example, we would see that the Earth pulls it towards the left, and the Moon pulls it to the right, and the resultant force is just right to keep the spacecraft in orbit between the Earth and the Moon.

---

[73] Euler, L., *Recherches sur le mouvement des corps célestes en général* (Studies on the movement of celestial bodies in general). (E112) Presented to the Berlin Academy on 8 June 1747. Published in Mémoires de l'Académie des sciences de Berlin 3, 1749, pp. 93-143

[74] Euler, L., *Theoria motus lunae exhibens omnes eius inaequalitates* (Theory of the motion of the moon which exhibits all its irregularities). Published as a book in 1753.

[75] Fellmann, E.A., *Leonhard Euler*. Birkhäuser Verlag, Switzerland (2007).

[76] Euler, L., *Recherches sur la question des inegalites du mouvement de Saturne et de Jupiter* (Investigations on the question of the differences of the movement of Saturn and Jupiter), Presented to the St. Petersburg Academy on 20 October 1768. Originally published as a book in 1772. The Euler Archive.

[77] Euler, L., *Considerationes de motu corporum coelestium* (Considerations on the motion of celestial bodies) E304. Read to the Berlin Academy on April 22, 1762.

Euler's pioneering effort to solve the three-body problem was the basis upon which other mathematicians, including Lagrange, built their own methods. The linear solutions to the equation of the fifth degree are called "Lagrange solutions" without mentioning Euler's contribution. However, Euler was the first to derive a theory of perturbations in 1762. By iteration, he determined, for the first time, the perturbation of the elements of the elliptical paths, and then applied this method to determine the motion of three mutually attracting bodies.[78]

In 1763, Euler considered the three bodies moving along a straight line.[79] He expanded the analysis and, two years later, he treated the collinear case in detail, presenting the differential equations of the *n*-body problem with ten integrals.[80] Regarding the solution of the three-body problem, Euler concluded with this statement:

> Mais, des qu'il est question de trois corps, dont le mouvement est déterminé par 9 équations, les sept équations intégrales que je viens de trouver ne suffisent plus pour en tirer une solution parfaite ; il faudrait encore au moins en découvrir deux nouvelles, auxquelles on n'a pu encore parvenir, malgré tous les soins que les plus grands Géomètres se sont donnés. La méthode dont je me suis servi ici, en cherchant certaines combinaisons entre les équations principales détaillées dans le § 35, qui conduisent a quelque équation intégrable, semble entièrement épuisée, & il faudra fans doute chercher une route tout à faire nouvelle. Dans l'état ou l'Analyse se trouve, il semble même impossible de dire si l'on en est encore fort éloigné ou non ; mais il est bien certain que, des qu'on sera arrivé à ce point, l'Analyse en retirera de beaucoup plus grands avantages, que l'Astronomie ne faudrait s'en promettre, à cause de la grande complication dont tous les éléments seront entrelaces selon toute apparence, de forte que pour la pratique on ne pourra presque en espérer aucun secours.

In 1770, during the most productive decade of his career despite his rapidly progressing blindness, Euler calculated the elliptical orbit of the comet observed in 1769,[81] and he won the 1770 Paris Academy prize for his *Théorie de la Lune*, which gave an unsatisfactory solution to the three-body problem. Euler shared the prize with his son Albrecht, who most likely simply took notes dictated to him by his blind father.

---

[78] J.J. Burckhardt, *Leonhard Euler 1707-1783*, Mathematics Magazine 56 (1983), 262-272. As it appeared in Sherlock Holmes in Babylon, Eds. M. Anderson, V. Katz, and R. Wilson, The Mathematical Association of America (2004).

[79] Euler, L., *De motu rectilineo trium corporum se mutuo attrahentium* (On the rectilinear motion of three bodies mutually attracted to each other). (E327). Presented to Petersburg Academy on 21 December 1763.

[80] Euler, L., *Considerations sur le probleme des trois corps*. (E400). Read to the Berlin Academy on 4 December 1765. Published in Memoires de l'Academie des sciences de Berlin 19, 1770, pp. 194-220

[81] Euler, L., *Recherches et calculs sur la vraie orbite elliptique de le comete de l'an 1769 et son tems periodique* (Research and calculations on the true elliptical orbit of the comet of the year 1769 and its periodic time). (E389). Presented to the St. Petersburg Academy on 10 September 1770.

The topic proposed by the Royal Academy of Sciences of Paris for the year 1772 was a follow-up to that which was issued two years earlier:

> De perfectionner les méthodes sur lesquelles est fondée la théorie de la Lune, de fixer par ce moyen celles des équations de ce Satellite, qui sont encore incertaines, & d'examiner en particulier si l'on peut rendre raison, par cette théorie de l'équation séculaire du mouvement moyen de la Lune.[82]

Once again Euler, now sixty-five, submitted an anonymous memoir: *Nouvelles Recherches sur le vrai mouvement de la Lune. Où l'on détermine toutes les inégalités auxquelles il est assujetti*. It carried the epitaph: "Hic labor extremus, longarum hæc viarum, hinc jam digissî, vestris appellimus oris." This is a passage from Virgil's *Aeneas*.[83]

Euler proposed the use of a rotating or synodic coordinate system for the restricted problem of the Sun-Earth-Moon system, and he concentrated on the true motion of the Moon, determining all inequalities to which it is subject. The younger Lagrange also focused on the Earth-Sun-Moon system, in effect solving the *restricted* three-body problem. Lagrange submitted *Essai sur le problème des trois corps*, aiming for a general solution of the three-body problem. He clarified that his method of solution was different from all those which had been given before. "It consists of not using in the determination of the orbit of each body any other elements other than the distances between the three bodies, that is to say, the triangle formed by these bodies at each instant." Then Lagrange stated that "The second Chapter examines how and in which cases the three bodies can move so that their distances shall always be constant, or shall at least keep constant ratios between themselves. I find that these conditions can only be satisfied in two cases: one, when the three bodies are arranged in one straight line, and the other, when they form an equilateral triangle; then each of the three bodies describes around the two others circles or conic sections, as if there were only two bodies."[84]

He demonstrated two special solutions: the collinear and the equilateral, for any three masses, with circular orbits. In this manner, Lagrange found two constant-pattern configurations for the three bodies where specific regions in the space surrounding the two large bodies like the Sun-Earth system within which a smaller mass can orbit naturally while maintaining the same position with respect to the other two. In 1762, Euler had already deduced the existence of two of those locations.

The combined analytical solution of Euler and Lagrange showed that there exist unique regions in space, which are five in every two-body system, which are equilibrium points for the system. Because Lagrange's solution was more comprehensive, those equilibrium points in space are now called Lagrangian points, and are identified as L1, L2, L3, L4 and L5. The 1772 prize of the Paris Academy was awarded to both Euler and Lagrange.

---

[82] Recueil des pièces qui ont remporté les prix de l'Académie royale des sciences. Tome Neuvieme. Qui contiens les Pièces de 1764, 1765, 1766, 1770 & 1772.

[83] Translation: "This was my last disaster, this is the termination of my long tedious voyage," Aeneas said at the death of his father, referring to it as the last trial at the end of his long wanderings. (3.714)

[84] Lagrange, J.-L., *Essai d'une Nouvelle Méthode pour résoudre le problème des trois corps*. (1772), pp. 229-230.

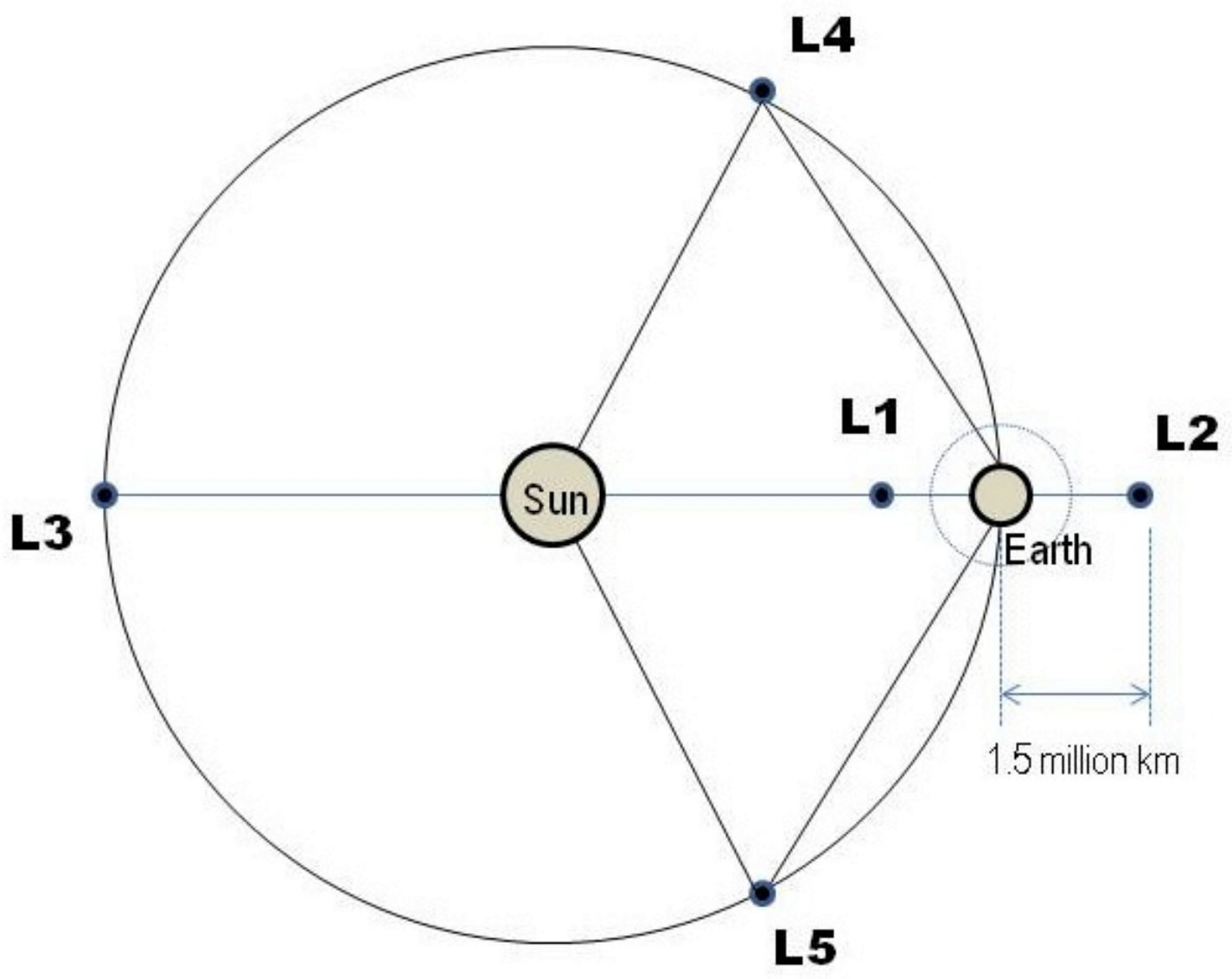

The five equilibrium points in a two-body system with one body (Sun) far more massive than the other (Earth). The L3–L5 locations appear to share the secondary body's orbit, however they are situated slightly outside it.

These five positions in an orbital configuration indicate where a small body, affected only by gravity, can theoretically be stationary relative to two larger bodies. In other words, these are the positions where the combined gravitational pull of the two large masses provides precisely the centripetal force required to rotate with them. One example of such system is that comprised of a spacecraft positioned with respect to the Earth and Moon.

Of the five equilibrium points, two are dynamically stable and three are unstable for small masses. This means that tiny departures from equilibrium will grow exponentially with time. For the Earth-Sun system, for example, a spacecraft orbiting at either of these two points will wander off after a few months, unless orbit corrections are made to keep the spacecraft in its position.[85] Note in the figure the unstable points (L1, L2 and L3) lie along the line connecting the two large masses. The stable points (L4 and L5) form the apex of two equilateral triangles that have the large masses at their vertices. Those triangular points L4 and L5 are called *libration* points because a body there will librate[86] instead of escaping, provided that mass ratio of the primaries is smaller than 0.0385. For example, for the Earth-Moon system, the mass ratio is 0.012. At the beginning of the twentieth century, astronomers discovered the Jupiter Trojan asteroids at the L4 and L5 points of the Sun–Jupiter system.

In order to model the Moon-Earth-Sun system, Euler introduced the orthogonal rotating coordinate system shown below. The schematic represents the ecliptic plane, the plane in which the Earth's orbit around the Sun lies or, more precisely, the plane in which the centre of gravity of the Earth-Moon system (its barycentre) orbits the Sun. Euler assumed the center of the stationary Sun at O (left side), and at a given time the Center of the Earth at ♀, with the right ♈ marking the location of the spring equinox. Either then out of this plan, ☾ the center of the Moon, where he lowered on this plane of the ecliptic perpendicular ☾V, of strong right ♀V marks the true longitude of the Moon, & ☾♀V the latitude angle, assumed to the north. Euler then shoot in the plane of the ecliptic the right line ♀M to indicate the mean longitude of the Moon, on which he lowered finally V perpendicular to VL, to have the three orthogonal coordinates L, LV, & V☾. This configuration is somewhat different from that which he used in the previous memoir; but Euler believed that in this manner he'd get the advantage that his formulas become there simpler and their calculation less awkward.

[85] For a full derivation of the Lagrangian points, refer to Neil J. Cornish analysis published on the web: www.physics.montana.edu/faculty/cornish/lagrange.html
[86] Libration is an oscillatory motion around an equilibrium point.

In a previous essay, Euler had defined the Cartesian coordinate system with respect to the ecliptic as ♀L = *a*(1+ *x*), LV = *ay*, and ☾V = *az*, where $a$ represented the average distance from the Moon to the Earth, and 1 expressed the average distance from the Earth to the Sun. For the 1772 memoir, he reduced the coordinates ♀L = 1+ *x*, LV = *y*, and ☾V = *z,* and kept the average distance from the Earth to the Sun expressed as 1/*a*, and from the known parallax of the Sun Euler determined 1 / *a* = 390.

This reference system rotates around Earth with the mean motion of the Moon marked by its true longitude ♁V.

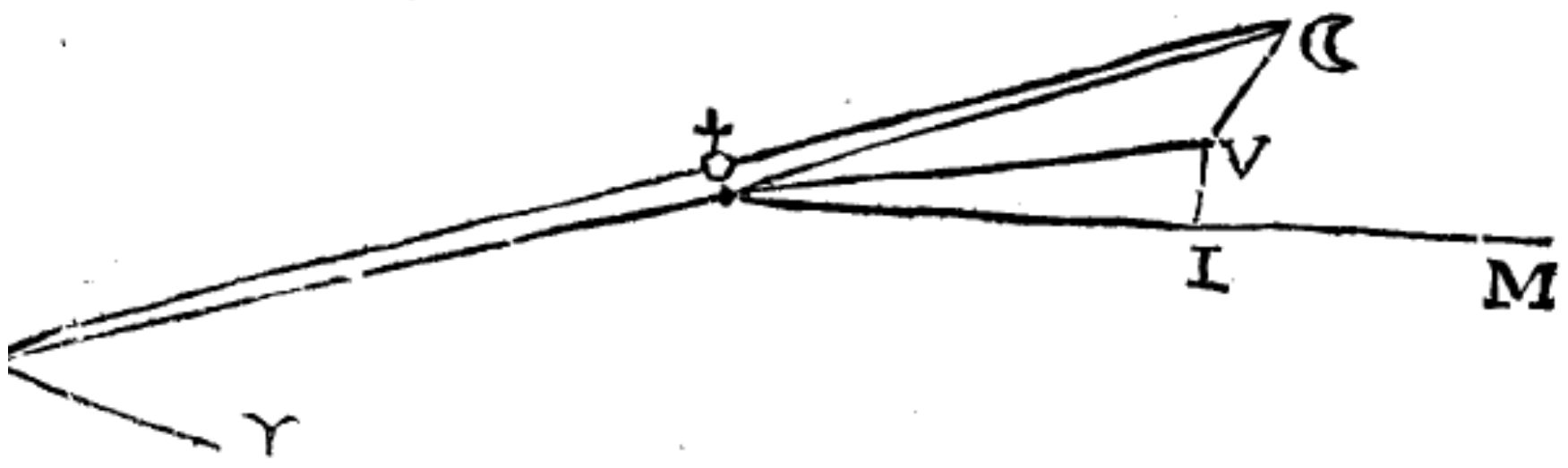

Coordinate system as depicted in Euler's essay of 1772.

The lines connecting the Moon to the Earth and to the Sum may have been Euler's attempt to indicate that the orbit of the Moon is inclined 5.145396° with regard to the ecliptic. Again, let us keep in mind that Euler was blind, and the sketches in his memoir were the interpretations of his assistants derived from his verbal explanations. In any case, Euler derived the second order partial differential equations for the three coordinates, each equation involving up to thirty-two terms:

$$\frac{d^2x}{dt^2} - \frac{2(m+1)dy}{dt} - 3\lambda x - \tfrac{3}{2}\cos 2p - \tfrac{3}{2}x\cos 2p + \ldots = 0$$

$$\frac{d^2y}{dt^2} + \frac{2(m+1)dx}{dt} + \tfrac{3}{2}\sin 2p + \tfrac{3}{2}x\sin 2p + \ldots = 0$$

$$\frac{d^2z}{dt^2} + (\lambda + 1)z - 3\lambda xz + 6\lambda x^2 z - \tfrac{3}{2}\lambda z(y^2 + z^2) - \cdots = 0$$

After solving his equations and comparing his results with observational data, Euler must have realized that his analysis did not, in effect, provide the required solution to the dynamics of the Moon-Earth-Sun system. He ended his memoir with these words: "Par-là, j'espère avoir pleinement satisfait aux vues de l'illustre Académie Royale des Sciences, ayant entièrement développé & fixé toutes les inégalités, auxquelles le mouvement de ce satellite est assujetti, sans en excepter aucune, qui pourrait influer sur le lieu de la Lune pour plus de dix secondes, & ayant enfin fait voir que, comme aucune de ces inégalités ne saurait produire une équation séculaire, dans le moyen mouvement de la Lune, on n'en pourra non plus rendre raison par la seule attraction du Soleil & de tout autre Corps céleste; de forte qu'il ne reste plus aucun doute que cette équation séculaire, qu'on observe, ne soit l'effet de la résistance du milieu, dans lequel les planètes se meuvent." Indeed, a higher degree of precision could not be obtained due to the action of the Sun and all the planets to which the Moon is subject.

Assisted by his son Johann Albrecht Euler, and his closest associates Kraftt and Lewell, Euler published his last lunar tables in 1772.

In 1836, German mathematician Carl Gustav Jacob Jacobi used Euler's rotating system to derive the integral of the motion named after him. The Jacobi integral is the constant of motion in the co-rotating frame of the restricted three-body problem, given by

$$C \equiv -2H = x^2 + y^2 + \frac{2(1-\mu)}{r_1} + \frac{2\mu}{r_2} - \dot{x}^2 - \dot{y}^2,$$

where $H$ is the Hamiltonian, $x$ and $y$ are the coordinates of the small mass $\mu$, and $r_1$ and $r_2$ are the distances from the two central masses.

It would take over a hundred years after the Euler-Lagrange prize, before French mathematician Henri Poincaré discovered that the full three-body problem is chaotic and cannot be solved in close form. And not until the last half of the twentieth century mathematicians and astronomers identified the non-irregularity of the planetary motion in the Solar System by using numerical models of stability over millions of years. Poincaré had established the non-integrability of the three-body problem in the late nineteenth century, however, not until recently we learned that the motion of the planets is chaotic, which requires accurate predictions of their trajectories over very long periods of time that correspond to the age of our Solar System.

Although there are many instances of three-body configurations in the Solar System that can be accurately described by the restricted three-body problem, such as when we deal with spacecraft moving between two large masses, for cases where it is not possible, then we use numerical methods to seek approximated solutions. Euler would agree.

Moreover, Euler would be thrilled to learn that, in the twentieth century, the space exploration program presented a host of challenging three-body problems that he helped resolve. For example, the initial problem was that of solving the motion of artificial satellites, followed by the task of determining spacecraft orbital maneuvers required for missions to the Moon. Today, two important artificial satellites orbit right at the two points in space that Euler predicted: the Solar and Heliospheric Observatory (SOHO) orbits around the first Lagrangian point.

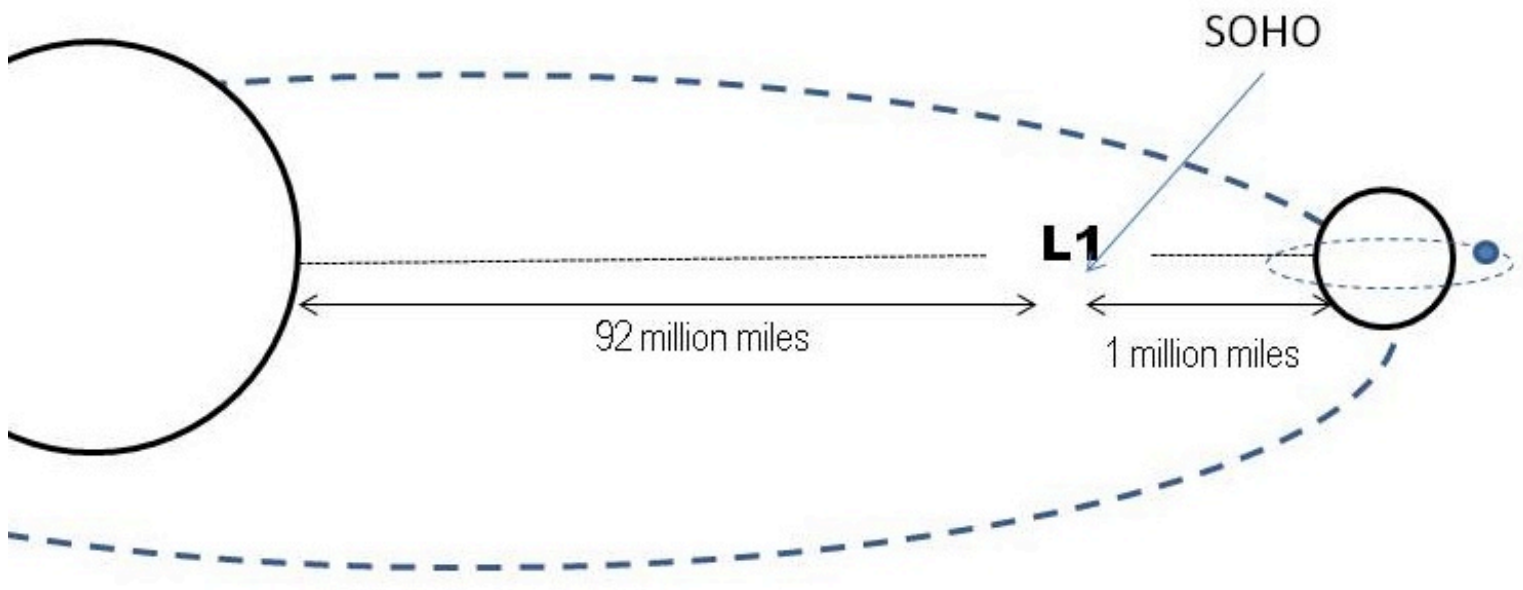

Launched in late 1995, SOHO is an international project of collaboration between ESA and NASA to study the Sun from its deep core to the outer corona and the solar wind. SOHO moves around the Sun in step with the Earth, by slowly orbiting around L1 where the combined gravity of the Earth and Sun keep SOHO in an orbit locked to the Earth-Sun line. L1 is a good position from which to monitor the Sun since the constant stream of particles from the Sun, the solar wind, reaches L1 about an hour before reaching Earth. Sunspots research is also part of SOHO's mission.

The Wilkinson Microwave Anisotropy Probe (WMAP) is positioned at L2. WMAP is a spacecraft that measures differences in the temperature of the Big Bang's remnant radiant heat—the Cosmic Microwave Background Radiation—across the full sky. As we noted above, the L1 and L2 points are unstable on a time scale of approximately twenty-three days. This requires that spacecraft parked at these positions undergo regular course and attitude corrections.

The future James Webb Space Telescope (JWST) will be put in an orbit around the L2 point. The reason is that the space telescope will observe primarily the infrared light from faint and very distant celestial objects. Because all objects, including the telescope, also emit infrared light, it is necessary that the telescope and its instruments be kept at very cold temperatures to avoid flooding the very faint astronomical signals with radiation from the telescope proper. Moreover, the JWST telescope will have a large shield to block the light from the Sun, Earth, and Moon, which otherwise would heat up the telescope, and interfere with the observations. That is why JWST must be in an orbit around L2 where all three bodies are in the same direction.

The space exploration program puts in perspective the complexity of the three-body problem. For analysis of space missions, such a spacecraft moving between planets, or from Earth to the Moon, we can apply the restricted three-body problem because only the small mass spacecraft is affected by the gravitational forces from two large bodies, and neither the Earth nor the Moon feel the influence of the spacecraft.

However, the complex space missions through the farther regions of the Solar System require that spacecraft trajectories be designed by decoupling an *n*-body system into several three-body systems. One interesting example is a mission to Jupiter. The spacecraft's trajectory orbit is designed by decoupling the Jovian moon *n*-body system into several three-body systems. The Jupiter-Ganymede-Europa-spacecraft four-body system can be approximated as two coupled planar three-body systems. In such approach, the two adjacent moons competing for control of the spacecraft are modeled as two nested three-body systems, e.g., Jupiter-Ganymede-spacecraft and Jupiter-Europa-spacecraft. Then, the desired orbit is designed with the two planar three-body systems. The initial solution is then refined to obtain a trajectory in a more accurate four-body model.[87]

## 6. Optics and Telescopes

Is it surprising that a blind mathematician would have wanted to improve the instruments required to observe the stars? Euler's work helped the progress of mathematical optics that led to the improvement of the telescope. Euler interest in optics began before he was thirty and remained almost to his death. He provided much of the needed mathematical analysis to understand light dispersion and chromatic aberration. Euler treated atmospheric refraction of light analytically both for celestial and terrestrial objects. With his mathematical optics, Euler also deduced the law of refraction.

In 1748, Euler wrote "On the perfection of objective lenses of telescopes." Searching for a way to correct chromatic aberration, Euler studied dispersion, which Newton had declared unattainable. Yet, Euler found the formula that the ratio of the natural logarithms of the indices of refraction must remain constant for different colors. This investigation provided the needed evidence that induced English optician John Dollond to construct his achromatic lenses. Euler calculated a large number of composite telescopes and microscopes, thereby increasing the possibilities in the construction of these instruments.

In 1769, Euler published *Dioptrica*, a three-volume textbook in which he collected his research in optics. He obtained the familiar formulae of elementary optics and stated a theory he believed could bring optical instruments to the highest degree of perfection. Euler studied the causes of failure in the poor quality of the lenses and attempted to determine errors in the laws of diffraction, which were determined experimentally. In 1774, the year after Swiss mathematician Nicolas Fuss arrived in St. Petersburg to work as assistant to Euler, he published a book on the construction of telescopes and microscopes, which contained a summary of Euler's contributions.

More than two centuries later, the University of Geneva in Switzerland named their modern 1.2-meter reflector telescope, the Euler Telescope in his honor. Built at the Geneva Observatory, located at La Silla (Chile), the Euler telescope is used in conjunction with the Coralie spectrograph to conduct high-precision radial velocity measurements principally to search for large extrasolar planets in the southern celestial hemisphere. Its first success was the discovery of a planet in orbit around Gliese 86, a K-type main-sequence star approximately 35 light-years away in the constellation of Eridanus. The planet (identified as Gliese 86b to distinguished from its parent star) orbits very close to the star, completing an orbit in 15.78 days.

To date, nearly two thousand planets have been confirmed to be orbiting other stars in our galaxy. Euler would be pleased to know, for it would confirm his belief. In 1760, Euler wrote: "Chaque étoile fixe semble être destinée pour échauffer & éclairer un certain nombre de corps opaques, semblables à notre terre, & habités aussi sans doute, lesquels se trouvent dans son voisinage, mais que nous ne voyons point à cause de leur prodigieux éloignement."[88]

[87] Gomes, et al., AAS 01-301, 2001

## 7. Euler's Legacy

Leonhard Euler lived in an era of great scientific development when the combination of theoretical models and observations led to a better understanding of the universe. In the eighteenth century, new instruments allowed the cataloging of stars with increasing precision and led to the discovery of new phenomena such as aberration, nutation, nebulae, and redefining the Solar System.

With better telescopes astronomers were able to look deeper into the night sky. In 1781, astronomer William Herschel discovered a new bright object, what he thought was a new comet, and named it Georgium Sidus in honor of his patron, King George III. However, the amateur German astronomer, Wilhelm Olbers had developed a method to calculate the orbits of comets, and determined that Herschel's object did not follow the same kind of orbit as comets. Olbers suggested that it was more likely to be a planet and not a comet. He was correct.

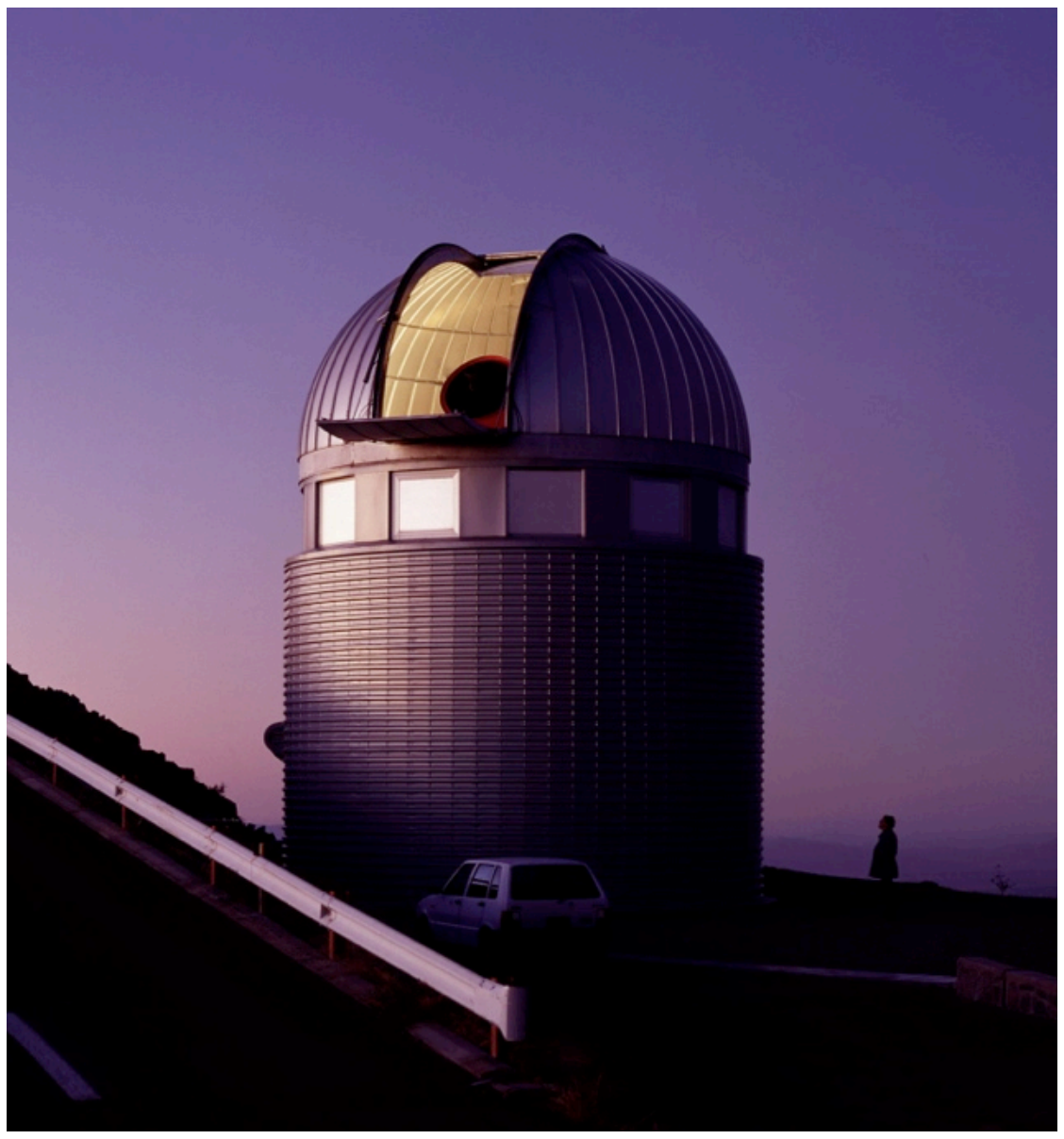

The dome of the Swiss Leonhard Euler Telescope at La Silla. Credit: ESO/H.Zodet

Euler's death is befitting a great master mathematician astronomer. On the 7th (18) of September 1783, the 76-year-old Euler was engaged in his mathematical research. After dinner he talked with his assistant Fuss about the discovery of a *seventh planet,* the new planet discovered by Herschel and that would be named Uranus.[89] Lexell had also visited Euler that day,[90] and the three talked about the analysis required to calculate its orbit. A little while later, while playing with his grandson, all of a sudden the pipe Euler was smoking slipped from his hand, and moments later with the words *Ich sterbe*! (I am dying!), lost consciousness;[3] Euler had suffered a stroke and died a few hours after. On that day, Euler "*il cessa de calculer et de vivre*."[91]

[88] Euler, L., *Lettres à une Princesse d'Allemagne sur divers sujets de Physique & de Philosophie.* A Mietau et Leipsic, 1770. Letter LIX, p. 247.

To honor Euler, astronomers have named a lunar crater and an asteroid after him. *Crater Euler* is a lunar impact crater 27 kilometers in diameter and about 2.5 kilometers deep. It is located in the southern half of the Mare Imbrium, with coordinates 23.3°N, 29.2°W. The most notable nearby feature near crater Euler is Mons Vinogradov to the west-southwest. There is a cluster of low ridges to the southwest, and this formation includes the small Natasha crater and the tiny Jehan crater.

On 29 August 1973, Russian astronomer Tamara Mikhailovna Smirnova,[92] discovered an object in the main asteroid belt that was named Asteroid *2002 Euler*. It has a diameter of 17.44 km. *2002 Euler* was at it closest approach to Earth (1.399 AU) in February 2013. The closest approach happens as a heavenly body has an orbit that somewhere is closest to Earth. Most asteroids and comets have an elliptical orbit in relation to a star (i.e. the Sun) or planet, and thus will vary in distance to Earth.

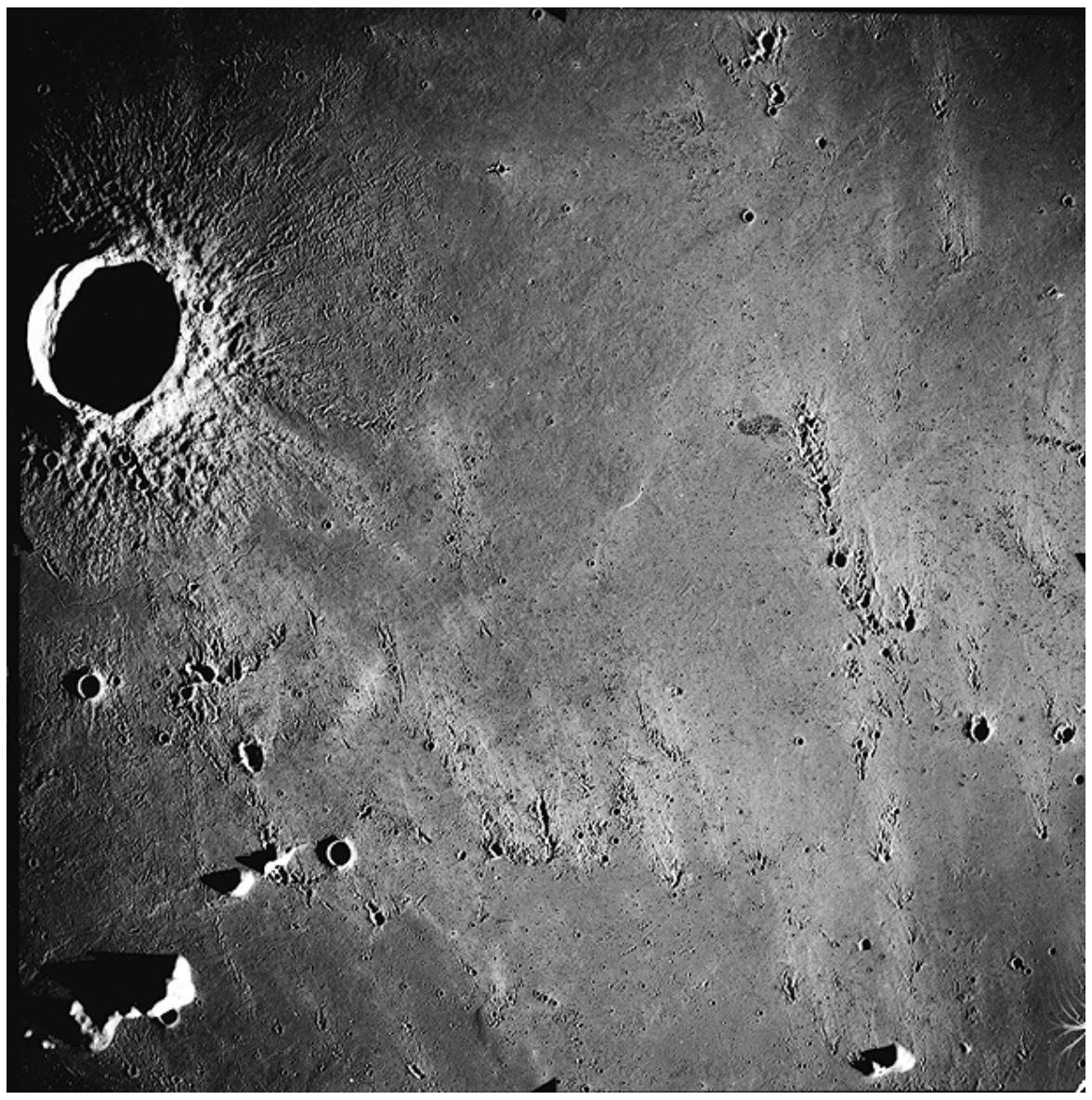

Euler crater near Copernicus secondary crater chains on the Moon. The image was taken from a distance of 114 km. NASA image: Apollo 17, AS17-2291. Credit: NASA.

[89] William Herschel discovered Uranus in 1781. While measuring the directions and brightness of stars, Herschel found a fuzzy spot that moved among the stars. This was Uranus, the first planet that was not known to the ancients

[90] Fellmann, E.A., *Leonhard Euler*, Birkhäuser Verlag, Switzerland (2007).

[91] Marquise de Condorcet, *Eloge de M. Euler*. In *Lettres de M. Euler a une Princesse d'Allemagne, sur différentes questions de physique et de philosophie*. Nouvelle Edition, Paris 1787. p. xliii.

[92] The asteroid 5540 Smirnova was named in her honor.

In addition to developing new mathematics and solving many problems in astronomy, Euler the blind mathematician was a great teacher. Not in the conventional manner since his work did not involve university teaching. Euler's career developed at the Imperial Academy of Sciences in St. Petersburg, and at the Royal Society of Sciences in Berlin. However, Euler taught many people, starting with his own son Johann Albrecht, who became a mathematician himself, and other researchers at the Academies. Euler mentored eight of sixteen professors belonging to the Russian Academy of St. Petersburg; all of them were well known and decorated with academic honours.

Euler's teaching became legendary through his "Letters to a German Princess on different topics of physics and philosophy." This is a three-volume book containing more than two hundred letters he wrote between 1760 and 1762 to instruct a fifteen-year old girl. Euler taught her many topics pertaining to physics and astronomy, including the notion of gravity and the laws that govern the motion of celestial bodies.[93] This book was translated in several languages and became widely read for decades after his death.

In the Eulogy, N. Fuss wrote: "The Academy has eight mathematicians who have enjoyed the instructions of Mr. Euler ... J. A. [Johann Albrecht] Euler, Kotelnikov, Rumovsky, Krafft, Lexell, Inokhodtsov, Golovine and I..." These individuals—some were astronomers—helped to carry on Euler's great legacy in many ways. Fuss served assistant to Euler from 1773–1783, and remained in Russia until his death in 1826. He contributed to spherical trigonometry, differential equations, the optics of microscopes and telescopes, differential geometry, and actuarial science.

French mathematician the marquis de Condorcet wrote in his eulogy, "All the noted mathematicians of the present day are his pupils: there is no one of them who has not formed himself by the study of his works, who has not received from him the formulas, the method which he employs; who is not directed and supported by the genius of Euler in his discoveries."

I disagree with Condorcet, because one does not have to be a *noted mathematician* to have been instructed by Euler, the universal scholar. Engineers, physicists, and astronomers, we owe much of what we know about analysis to Euler, and thus, we all are Euler's pupils.

## Acknowledgements

Many thanks to Dominic Klyve (Central Washington University), Lee Stemkoski (Adelphi University), and Erik Tou (Pacific Lutheran University), for making available online Euler's extensive publications at the **Euler's Archive**, and to the Mathematical Association of America for hosting it at http://eulerarchive.maa.org/.

---

[93] Musielak, D., *Euler and the German Princess*, Freiburg, Germany, 28 June 2014. This article provides insights regarding the life of the princess and reviews some of the material in the letters.

# Euler Timeline

| | | |
|---|---|---|
| Basel & Riehen | 1707 | Euler is born in Basel, Switzerland |
| | 1720 | Enters University of Basel |
| | 1723 | Completes master degree in philosophy |
| St. Petersburg, Russia | 1727 | Begins tenure at St. Petersburg Academy. Newton dies early the same year. Euler Begins daily astronomical observations at St Petersburg Observatory; activity would last 10 years |
| | 1730 | Appointed professor of physics at St. Petersburg Academy |
| | 1732 | Publishes *Solution to problems of astronomy*, an elementary method to calculate stellar coordinates |
| | 1734 | Publishes first volume of his monumental book *Mechanics*. Euler outlined a program of studies embracing every branch of science, involving a systematic application of analysis. He laid the foundations of analytical mechanics. |
| | 1735 | Publishes *On the motion of planets and orbits,* where he introduces an iterative method to solve Kepler's problem. |
| | 1736 | Publishes the second volume of *Mechanica* |
| | 1738 | Looses vision in his right eye |
| | 1740 | Wins Paris Academy **Prize** for work explaining the causes of the Tides. He shares prize with Daniel Bernoulli |
| Berlin | 1741 | Begins his tenure at the Berlin Academy of Sciences |
| | 1744 | Publishes *Theory of the Planets and Comets*, his first astronomy book |
| | 1747 | Formalizes restricted three-body problem |
| | 1748 | Wins Paris Academy **Prize** for research on Jupiter, Saturn inequalities |
| | 1752 | Wins Paris Academy **Prize** for "*Sur les inégalités de Jupiter*." |
| | 1753 | Publishes an 18-chapter book titled *Theory of the motion of the moon which exhibits all its irregularities* |
| | 1756 | Wins Paris Academy **Prize** for his work on "*Sur les inégalités de la Terre*." |
| St. Petersburg (Russia) | 1766 | Returns to Russian Royal Academy of Sciences |
| | 1768 | Wins Paris Academy **Prize** for his second Moon theory, coauthored with his son Johann Albrecht |
| | 1770 | Wins Paris Academy **Prize** for extension of his second Moon theory, sharing honor with his son Johann Albrecht |
| | 1771 | Becomes completely blind |
| | 1772 | Wins the Paris Academy **Prize** for the Moon theory or three-body problem. Lagrange shares the prize for his refined solution. |
| | 1772 | Publishes *Theoria motum Lunae* |
| | 1777 | Euler's thoughts concerning the Earth's motional inequalities caused by Venus' action, accompanied with a table of the corrections of Earth's position |
| | 1783 | Euler ceases to calculate. |